\input ./style/arxiv-general.cfg
\documentclass[MSNbibl,number,citesort,dvips]{arxbj}
\makeatletter
   \@ifpackageloaded{graphicx}{}{\usepackage{graphicx}}
\makeatother

\aid{0}
\volume{21}
\issue{4}
\pubyear{2015}
\firstpage{1984}
\lastpage{2023}
\doi{10.3150/14-BEJ633} 

\makeatletter
\mathchardef\varPi="0105
\mathchardef\varPsi="0109
\newcommand{\mh}{\mathfrak{h}}
\newcommand{\rrvert}{\vert}
\newcommand{\rrVert}{\Vert}
\newcommand{\llvert}{\vert}
\newcommand{\llVert}{\Vert}
\def\cal{\mathcal}
\newtheorem{lemma}{Lemma}
\newtheorem{theorem}{Theorem}
\newtheorem{proposition}{Proposition}
\newremark{example}{Example}
\newremark{remark}{Remark}
\newproclaim{definition}{Definition}
\newtheorem{corollary}{Corollary}

\renewcommand{\kappa}{\varkappa}
\newcommand{\rd}{\mathrm{d}}
\newcommand{\rh}{\mathrm{h}}

\newcommand{\cA}{{\cal A}}
\newcommand{\cB}{{\cal B}}
\newcommand{\cH}{{\cal H}}
\newcommand{\cI}{{\cal I}}
\newcommand{\cN}{{\cal N}}
\newcommand{\cP}{{\cal P}}
\newcommand{\cR}{{\cal R}}
\newcommand{\cU}{{\cal U}}

\newcommand{\bE}{\mathbb E}
\newcommand{\bF}{\mathbb F}
\newcommand{\bH}{\mathbb H}
\newcommand{\bL}{{\mathbb L}}
\newcommand{\bN}{{\mathbb N}}
\newcommand{\bP}{{\mathbb P}}
\newcommand{\bR}{{\mathbb R}}
\newcommand{\bW}{{\mathbb W}}
\newcommand{\mA}{\mathfrak{A}}
\newcommand{\mB}{\mathfrak{B}}
\newcommand{\mF}{\mathfrak{F}}
\newcommand{\mH}{\mathfrak{H}}
\newcommand{\mP}{\mathfrak{P}}
\newcommand{\mz}{\mathfrak{z}}
\newcommand{\mR}{\mathfrak{R}}

\makeatother

\begin{document}
\begin{frontmatter}

\title{Pointwise adaptive estimation of a~multivariate density under independence~hypothesis}
\runtitle{Pointwise adaptive estimation}

\begin{aug}
\author[A]{\inits{G.}\fnms{Gilles}~\snm{Rebelles}\corref{}\ead[label=e1]{rebelles.gilles@neuf.fr}}
%
\address[A]{Institut de Math\'ematique de Marseille,
Aix-Marseille Universit\'e,
39, rue F. Joliot-Curie,
13453 Marseille, France. \printead{e1}}
\end{aug}

\received{\smonth{6} \syear{2013}}
\revised{\smonth{3} \syear{2014}}

%
\begin{abstract}
In this paper, we study the problem of pointwise estimation of a
multivariate density. We provide a data-driven selection rule from the
family of kernel estimators and derive for it a pointwise oracle
inequality. Using the latter bound, we show that the proposed estimator
is minimax and minimax adaptive over the scale of anisotropic Nikolskii
classes. It is important to emphasize that our estimation method
adjusts automatically to eventual independence structure of the
underlying density. This, in its turn, allows to reduce significantly
the influence of the dimension on the accuracy of estimation (curse of
dimensionality). The main technical tools used in our considerations
are pointwise uniform bounds of empirical processes developed recently
in Lepski [\textit{Math. Methods Statist.}
\textbf{22} (2013)
83--99].
\end{abstract}

%
\end{frontmatter}

\section{Introduction}
\label{sec:introduction}
Let $X_i=(X_{i,1},\ldots,X_{i,d}), i\in\bN^*$, be a sequence of $\bR
^d$-valued i.i.d. random vectors defined on a complete probability
space $ (\Omega,\mA,\textsf{P} )$ and having the density $f$
with respect to the Lebesgue measure. Furthermore, $\bP_{f}^{(n)}$
denotes the probability law of $X^{(n)}=(X_1,\ldots,X_n), n\in\bN^*$,
and $\bE_{f}^{(n)}$ is the mathematical expectation with respect to
$\bP
_{f}^{(n)}$.

Our goal is to estimate the density $f$ at a given point $x_0\in\bR^d$
using the observation $X^{(n)}=(X_1,\ldots,X_n)$, $n\in\bN^*$. As an
estimator, we mean any $X^{(n)}$-measurable mapping $\widehat{f}\dvtx\bR
^n\rightarrow\bR$ and the accuracy of an estimator is measured by the
\textit{pointwise risk}:
\[
\cR_n^{(q)} [\widehat{f},f ]:= \bigl(\bE_f^{(n)}
\bigl\llvert \widehat {f}(x_0)-f(x_0)\bigr\rrvert
^q \bigr)^{{1}/{q}}, \qquad q\geq1. %
\]

The discussion of traditional methods and a part of the vast literature
on the theory and application of the density estimation is given by
Devroye and Gy\"orfi \cite{devroye-gyorfi}, Silverman \cite{silverman}
and Scott~\cite{scott}.
We do not pretend here to provide with a detailed overview and mention
only the results which are relevant for considered problems.
The minimax and adaptive minimax multivariate density estimation with
$\bL_p$-loss on particular functional classes was studied in
Bretagnolle and Huber \cite{bretagnol}, Ibragimov and Khasminskii
\cite{ibragimov-khas1,ibragimov-khas2}, Devroye and Lugosi
\cite{devroye-lugosi1,devroye-lugosi2,devroye-lugosi3},
Efroimovich \cite{efroimovich1,efroimovich2}, Hasminskii and
Ibragimov \cite{hasminski-ibragimov}, Golubev \cite{golubev}, Donoho \textit{et
al.} \cite{dono2}, Kerkyacharian, Picard and Tribouley \cite{kerk1},
Gin\'e and Guillou \cite{Gine1}, Juditsky and Lambert-Lacroix \cite{judistky},
 Rigollet~\cite{rigollet}, Massart~\cite{massart}
(Chapter~7), Samarov and Tsybakov \cite{samarov-tsybakov}, Birg\'e
\cite{birge}, Mason \cite{mason}, Gin\'e and Nickl~\cite{Gine2},
Chac\'on
and Duong \cite{chaconduong} and Goldenshluger and
Lepski \cite{G-L:density-Lploss}. In Comte and Lacour~\cite{comte:deconvolution},
the pointwise setting was first considered in the context of
multidimensional deconvolution model. More recently, in Goldenshluger
and Lepski \cite{G-L:density-Lploss-bis}, adaptive minimax upper bounds
were proved for multivariate density estimation with $\bL_p$-risks on
anisotropic Nikolskii classes using a local (pointwise) procedure. The
use of Nikolskii classes allows to consider the estimation of
anisotropic and inhomogeneous densities; see Ibragimov and Khasminskii
\cite{ibragimov-khas2}, Goldenshluger and Lepski \cite{G-L:density-Lploss} and Lepski \cite{lepski:supnormlossdensityestimation}.

In this paper, we focus on the problem of the minimax and adaptive
minimax pointwise multivariate density estimation over the scale of
anisotropic Nikolskii classes.

\textit{Minimax estimation}. In the framework of the minimax
estimation, it is assumed that $f$ belongs to a certain set of
functions $\Sigma$, and then the accuracy of an estimator $\widehat{f}$
is measured by its \textit{maximal risk} over $\Sigma$:
%
\begin{equation}
\label{eq:risk} \cR_n^{(q)} [\widehat{f},\Sigma ]:=\sup
_{f\in\Sigma
} \bigl(\bE _f^{(n)}\bigl\llvert
\widehat{f}(x_0)-f(x_0)\bigr\rrvert ^q
\bigr)^{{1}/{q}},\qquad q\geq1.
\end{equation}
The objective here is to construct an estimator $\widehat{f}_*$ which
achieves the asymptotic of \textit{the minimax risk} (minimax rate of
convergence):
\[
\cR_n^{(q)} [\widehat{f}_*,\Sigma ]\asymp\inf
_{\widehat
{f}}\cR _n^{(q)} [\widehat{f},\Sigma ]:=
\varphi_n(\Sigma). %
\]
Here, infimum is taken over all possible estimators.

\textit{Smoothness assumption}. Let $\Sigma$ be either H\"older
classes $\bH(\beta,L)$ or $\bL_p$-Sobolev classes $\bW(\beta,p,L)$ of
univariate functions. Here, $\beta$ represents the smoothness of the
underlying density and $p$ is the index of the norm where the
smoothness is measured. Then
%
\begin{eqnarray}
\label{eq:univariateminimaxrates} \varphi_n \bigl(\bH(\beta,L) \bigr)&=&
n^{-{\beta}/{(2\beta+1)}},
\nonumber
\\[-8pt]
\\[-8pt]
\nonumber
\varphi_n \bigl(\bW(\beta,p,L) \bigr)&=&
n^{-{(\beta-1/p)}/{(2(\beta-1/p)+1)}},\qquad \beta>0, 1<p<
\infty.
\end{eqnarray}
These minimax rates can be obtained from the results developed by
Donoho and Low \cite{dono1}; see also Ibragimov and Khasminskii
\cite{ibragimov-khas1,ibragimov-khas2}, and Hasminskii and
Ibragimov \cite{hasminski-ibragimov}.

Let now $\Sigma=\bH_d(\beta,L)$ where $\bH_d(\beta,L)$ is an
anisotropic H\"older class determined by the smoothness parameter
$\beta
=(\beta_1,\ldots,\beta_d)$. In this case,
%
\begin{equation}
\label{eq:multivariateminimaxrates} \varphi_n \bigl(\bH_d(\beta,L)
\bigr)=n^{-{\overline{\beta
}}/{(2\overline{\beta}+1)}}, \qquad\overline{\beta}:= \Biggl[\sum
_{i=1}^d 1/\beta_i
\Biggr]^{-1}, \qquad\beta_i>0, i=\overline{1,d}.
\end{equation}
The latter result can be obtained from Kerkyacharian, Lepski and Picard
\cite{kerk2}, Proposition~1, in the framework of the Gaussian white
noise model. The similar minimax results will be established for
pointwise multivariate density estimation in Section~\ref{sec:minimaxresults}; see Theorems~\ref{theo:minimaxlowerbound1} and~\ref{theo:minimaxupperbound}.

It is important to emphasize that minimax rates depend heavily on the
dimension $d$. Let us briefly discuss how to reduce the influence of
the dimension on the accuracy of estimation (curse of dimensionality).
The approach which have been recently proposed in Lepski \cite{lepski:supnormlossdensityestimation} is to take into account the
eventual independence structure of the underlying density.

\textit{Structural assumption}. Note $\cI_d$ the set of all subsets
of $\{1,\ldots,d\}$ and $\mP$ the set of all partitions of $\{
1,\ldots,d\}$ completed by the empty set $\varnothing$. For all $I\in\cI_d$
and $x\in\bR^d$ note also $x_I=(x_i)_{i\in I}$, $\overline{I}=\{
1,\ldots,d\} \setminus I$, $\llvert I\rrvert =\operatorname{card}(I)$ and put
\[
f_I(x_I):=\int_{\bR^{\llvert \overline{I}\rrvert }}f(x)\,\mathrm{d}x_{\overline{I}}.
\]
Obviously, $f_I$ is the marginal density of $X_{1,I}$ and, to take into
account the independence structure of the density $f$, we consider the
following set:
\[
\mP(f):= \biggl\{\cP\in\mP\dvt f(x)=\prod_{I\in\cP}f_I(x_I),
\forall x\in \bR^d \biggr\}. %
\]

In this paper, we focus on the problem of pointwise multivariate
density estimation on anisotropic Nikolskii classes. In particular, we
will prove that the minimax rate on the class $N_{p,d}^*(\beta,L,\cP)$
(introduced in Lepski \cite{lepski:supnormlossdensityestimation}, see
the definition in Section~\ref{sec:nikolskiiclasses}) for fixed $\beta
\in(0,+\infty)^d$, $p\in[1,+\infty]^d$, $L\in(0,+\infty)^d$, $\cP
\in
\mP(f)$, are given by
\[
\varphi_n \bigl( N_{p,d}^*(\beta,L,\cP)
\bigr)=n^{-{r}/{(2r+1)}},\qquad
r:=\inf_{I\in\cP} \biggl[\frac{1-\sum_{i\in I}
{1}/{(\beta_i p_i)}}{\sum_{i\in I}{1}/{\beta_i}}
\biggr]. %
\]
If $d=1$, then the structural assumption does not exist, that means
formally $\cP=\overline{\varnothing}$, and we come to the rates given
in (\ref{eq:univariateminimaxrates}). Note that $N_{\infty,1}^*(\beta,L,\overline{\varnothing})$ coincides with the set of densities
belonging to $\bH(\beta,L)$ and that $N_{p,1}^*(\beta,L,\overline
{\varnothing})$ contains the set of densities belonging to $\bW(\beta,p,L)$.

If $d\geq2$, $p_i=\infty, i=\overline{1,d}$, and $\cP=\overline
{\varnothing}$ we find again the rates given in (\ref
{eq:multivariateminimaxrates}), and $N_{\infty,d}^*(\beta,L,\overline
{\varnothing})$ coincides with a set of densities belonging to $\bH
_d(\beta,L)$. Note however that if $\cP\neq\overline{\varnothing}$ the
latter rates can be essentially improved. Indeed, if, for instance,
$\beta
=(\bolds{\beta},\ldots,\bolds{\beta})$ and $\cP^*=
\{
\{1 \},\ldots, \{d \} \}$, then $r=\bolds
{\beta
}$ and
%
\begin{equation}
\label{comparison}
n^{-{{\beta}}/{(2\bolds{\beta}+d)}}=\varphi _n \bigl( \bH_d(
\beta,L) \bigr)\gg\varphi_n \bigl( N_{\infty,d}^*\bigl(\beta,L,
\cP ^*\bigr) \bigr)=n^{-{\bolds{\beta}}/{(2\bolds{\beta}+1)}}.
\end{equation}
Moreover, $\varphi_n ( N_{\infty,d}^*(\beta,L,\cP^*) )$
does not
depend on the dimension $d$.

We remark that minimax rates (accuracy of estimation) depend heavily on
the parameters $\beta, p$ and $\cP$. Their knowledge cannot be often
supposed in particular practice. It makes necessary to find an
estimator whose construction would be parameter's free.

\textit{Adaptive minimax estimation}. In the framework of the
adaptive minimax estimation the underlying density $f$ is supposed to
belong to the given scale of functional classes $ \{\Sigma
_{\alpha
}, \alpha\in\cA \}$. For instance, if $\Sigma_{\alpha}=\bH
(\beta,L)$, $\alpha=(\beta,L)$, or if $\Sigma_{\alpha}=\bW(\beta,p,L)$,
$\alpha=(\beta,p,L)$.

The first question arising in the framework of the adaptive approach
consists in the following: does there exists an estimator $\widehat
{f}_*$ such that
%
\begin{equation}
\label{eq:adaptivoptimal} \limsup_{n\rightarrow+\infty} \bigl\{\varphi_n^{-1}(
\alpha)\cR _n^{(q)} [\widehat{f}_*,\Sigma_{\alpha} ]
\bigr\} <+\infty \qquad\forall\alpha\in\cA,
\end{equation}
where $\varphi_n(\alpha)$ is the minimax rate of convergence over
$\Sigma_{\alpha}$.

As it was shown in Lepski \cite{lepski:pointwisembbadapt} for the
Gaussian white noise model, the answer of this question is negative if
$\Sigma_{\alpha}=\bH(\beta,L)$, $\alpha=(\beta,L)$. Brown and Low
\cite{brown} extended this result to the pointwise density estimation.
Further Butucea \cite{butucea:pointwisedensityadapt} extended the
results of Brown and Low \cite{brown} over the scale of $\bL_p$-Sobolev
classes $\bW(\beta,p,L)$. In Section~\ref{sec:optimality}, we will
prove that the answer is also negative for multivariate density
estimation at a given point over the scale of anisotropic Nikolskii
classes $N_{p,d}^*(\beta,L,\cP)$.\vspace*{1pt}

Thus, for problems in which (\ref{eq:adaptivoptimal}) does not hold we
need first to find a family of normalizations $\varPsi= \{\varPsi
_n(\Sigma_{\alpha}), \alpha\in\cA \}$ and an estimator
$\widehat
{f}_{\Psi}$ such that
%
\begin{equation}
\label{eq:adaptiv} \limsup_{n\rightarrow+\infty} \bigl\{\varPsi_n^{-1}(
\alpha)\cR _n^{(q)} [\widehat{f}_{\Psi},
\Sigma_{\alpha} ] \bigr\} <+\infty \qquad\forall\alpha\in\cA.
\end{equation}
Any family of normalizations satisfying (\ref{eq:adaptiv}) is called
admissible and the estimator $\widehat{f}_{\varPsi}$ is called
$\varPsi
$-adaptive. Next, we have to provide with the criterion of optimality
allowing to select ``the best'' admissible family of normalizations,
usually called adaptive rate of convergence. The first criterion was
proposed in Lepski \cite{lepski:pointwisembbadapt} and it was improved
later in Tsybakov \cite{tsybakov} and in Klutchnikoff~\cite{klutch}.

In particular, in Lepski \cite{lepski:pointwisembbadapt} and in Butucea
\cite{butucea:pointwisedensityadapt}, it was shown that the adaptive
rate of convergence for the considered problem is
\begin{eqnarray*}
\varPsi_n \bigl(\bH(\beta,L) \bigr)&=&\cases{ %
\displaystyle \biggl(
\frac{\ln(n)}{n} \biggr)^{{\beta}/{(2\beta+1)}}, &\quad
$\beta \in (0,\beta_{\mathrm{max}}),$
\vspace*{2pt}
\cr
\displaystyle\biggl(\frac{1}{n} \biggr)^{{\beta}/{(2\beta+1)}}, &\quad
$\beta=\beta_{\mathrm{max}},$}
\\
\varPsi_n \bigl(\bW(\beta,p,L) \bigr)&=&\cases{ %
 \displaystyle\biggl(
\frac{\ln(n)}{n} \biggr)^{{(\beta-1/p)}/{(2(\beta
-1/p)+1)}}, &\quad $\beta \in(0,\beta_{\mathrm{max}}),$
\vspace*{2pt}
\cr
\displaystyle\biggl(\frac{1}{n} \biggr)^{{(\beta-1/p)}/{(2(\beta-1/p)+1)}}, &$\beta =
\beta_{\mathrm{max}},$}
\end{eqnarray*}
with respect to the criterion in Lepski \cite{lepski:pointwisembbadapt}
and Tsybakov \cite{tsybakov}, respectively.
Here, $\beta_{\mathrm{max}}$ is an arbitrary positive number.

Later Klutchnikoff \cite{klutch} studied the pointwise adaptive minimax
estimation over anisotropic H\"older classes, in the Gaussian white
noise model. The consideration of anisotropic functional classes
required to develop a new criterion of optimality. Following this
criterion, Klutchnikoff~\cite{klutch} proved that the adaptive rate of
convergence is
\[
\varPsi_n \bigl(\bH_d(\beta,L) \bigr)=\cases{
 \displaystyle\biggl(\frac{\ln(n)}{n} \biggr)^{{\overline{\beta
}}{(2\overline{\beta
}+1)}}, &\quad $\displaystyle\beta\in\prod
_{i=1}^d\bigl(0,\beta_i^{(\mathrm{max})}
\bigr)^d,$ \vspace*{2pt}
\cr
\displaystyle\biggl(\frac{1}{n}
\biggr)^{{\overline{\beta
}^{(\mathrm{max})}}/{(2\overline
{\beta}^{(\mathrm{max})}+1)}}, &\quad $\beta=\beta^{(\mathrm{max})}.$} %
\]
Recently, Comte and Lacour \cite{comte:deconvolution} found a similar
form of admissible sequence for pointwise adaptive minimax estimation
in the deconvolution model.

In Section~\ref{sec:adaptiveminimaxresults}, we provide with minimax
adaptive estimator in pointwise multivariate density estimation over
the scale of anisotropic Nikolskii classes. We will take into account
not only the approximation properties of the underlying density but the
eventual independence structure as well. To analyze the accuracy of the
proposed estimator, we establish so-called pointwise oracle inequality
proved in Section~\ref{oracle-inequality-proof}.
We will also show that the adaptive rate of convergence is given by
\begin{eqnarray*}
\varPsi_n \bigl( N_{p,d}^*(\beta,L,\cP) \bigr)&=&\cases{
 \displaystyle\biggl(\frac{\ln(n)}{n} \biggr)^{{r}/{(2r+1)}},
 &\quad $0<r<r_{\mathrm{max}},$
\vspace*{2pt}
\cr
\displaystyle\biggl(\frac{1}{n} \biggr)^{{r}/{(2r+1)}},
&\quad $r=r_{\mathrm{max}}$,}\\
 r&:=&\inf_{I\in\cP} \biggl[
\frac{1-\sum_{i\in I}{1}/{(\beta_i
p_i)}}{\sum_{i\in I}{1}/{\beta_i}} \biggr]. %
\end{eqnarray*}
To assert the optimality of this family of normalizations, we
generalize the criterion proposed in Klutchnikoff \cite{klutch}; see
Section~\ref{sec:optimality}.

\textit{Organization of the paper}. In Section~\ref{selection-oracle-inequality}, we provide a measurable data-driven
selection rule based on bandwidth selection of kernel estimators and we
derive an oracle-type inequality for the selected estimator at a given
point. In Section~\ref{minimax-adaptive-estimation}, we treat the
complete problem of minimax and adaptive minimax pointwise multivariate
density estimation on a scale of anisotropic Nikolskii classes taking
into account the independence structure of the underlying density. In
Section~\ref{sec:discussion}, we briefly compare our local method with
the global one developed in Lepski \cite{lepski:supnormlossdensityestimation}. Proofs of all main results are
given in Section~\ref{proofs}. Proofs of technical lemmas are postponed
to the \hyperref[app]{Appendix}.

\section{Selection rule and pointwise oracle-type inequality}
\label{selection-oracle-inequality}
\subsection{Kernel estimators related to independence structure}
\label{sec:kernel}
Let $\mathbf{K}\dvtx\bR\rightarrow\bR$ be a fixed symmetric kernel
satisfying $\int{\mathbf{K}}=1$, $\operatorname{supp}(\mathbf{K})\subseteq
[-1/2,1/2]$, \mbox{$\llVert \mathbf{K}\rrVert _{\infty}<\infty$},
%
\begin{equation}
\label{eq:kernelasumption} \exists L_{\mathbf{K}}>0\dvt\qquad
 \bigl\llvert \mathbf{K}(x)-
\mathbf{K}(y)\bigr\rrvert \leq L_{\mathbf{K}}\llvert x-y\rrvert\qquad\forall x,y
\in\bR.
\end{equation}
For all $I\in\cI_d$, $h\in(0,1]^d$ and $x\in\bR^d$ put also
\begin{eqnarray*}
K^{(I)}(x_I)&:=&\prod_{i\in I}
\mathbf{K}(x_i),\qquad V_{h_I}:=\prod
_{i\in
I}h_i,\qquad K_{h_I}^{(I)}(x_I):=V_{h_I}^{-1}
\prod_{i\in I}\mathbf{K}(x_i/h_i);
\\
\widehat{f}_{h_I}^{(n)}(x_{0,I})&:=&n^{-1}
\sum_{i=1}^{n}K_{h_I}^{(I)}
(X_{i,I}-x_{0,I} ). %
\end{eqnarray*}
Then introduce the family of estimators
\[
\mF [ \mP ]:= \biggl\{\widehat{f}_{(h,\cP
)}^{(n)}(x_{0})=
\prod_{I\in\cP}\widehat{f}_{h_I}^{(n)}(x_{0,I}),
(h,\cP)\in (0,1]^d\times\mP \biggr\}. %
\]

Note first that $\widehat{f}_{(h,\overline{\varnothing
})}^{(n)}(x_{0})=\widehat{f}_{h}^{(n)}(x_{0})$ is the
Parzen--Rosenblatt estimator (see, e.g., Rosenblatt~\cite{rosenblatt},
Parzen \cite{parzen}) with kernel $K^{(\overline{\varnothing})}$ and
multibandwidth $h$.

Next, the introduction\vspace*{1pt} of the estimator $\widehat{f}_{(h,\cP
)}^{(n)}(x_{0})$ is based on the following simple observation. If there
exists $\cP\in\mP(f)$, the idea is to estimate separately each marginal
density corresponding to $I\in\cP$. Since the estimated density
possesses the product structure, we seek its estimator in the same form.

Below we propose a data driven selection from the family $\mF [
\mP
 ]$.

\subsection{Auxiliary estimators and extra parameters}
\label{oracle-notations-constants}
To define our selection rule, we need to introduce some notation and quantities.

\textit{Auxiliary estimators.}
For $I\in\cI_d$ and $h\in(0,1]^d$ put
\[
\widetilde{G}_{h_I}(x_{0,I}):=1\vee \Biggl[n^{-1}
\sum_{i=1}^{n}\bigl\llvert
K_{h_I}^{(I)} (X_{i,I}-x_{0,I} )\bigr
\rrvert \Biggr]. %
\]

Introduce for $I\in\cI_d$ and $h,\eta\in(0,1]^d$ auxiliary estimators
\[
\widehat{f}_{h_I,\eta_I}^{(n)}(x_{0,I}):=n^{-1}
\sum_{i=1}^{n}K_{h_I\vee
\eta_I}^{(I)}
(X_{i,I}-x_{0,I} ),\qquad h_I\vee
\eta_I:= (h_i\vee\eta_i
)_{i\in I}. %
\]
Note that the idea to use such auxiliary estimators, defined with the
multibandwidth $h\vee\eta$, appeared for the first time in
Kerkyacharian, Lepski and Picard \cite{kerk2}, in the framework of the
Gaussian white noise model.

We endow the set $\mP$ with the operation ``$\circ$'' introduced in
Lepski \cite{lepski:supnormlossdensityestimation}: for any $\cP,\cP
'\in
\mP$
\[
\cP\circ\cP':= \bigl\{I\cap I'\neq\varnothing, I\in
\cP, I'\in \cP' \bigr\}\in\mP. %
\]
Then we define for $h,\eta\in(0,1]^d$ and $\cP,\cP'\in\mP$
%
\begin{equation}
\label{eq:auxiliaryestimators} \widehat{f}_{(h,\cP),(\eta,\cP')}^{(n)}(x_{0}):=
\prod_{I\in\cP
\circ\cP
'}\widehat{f}_{h_{I},\eta_{I}}^{(n)}(x_{0,I}).
\end{equation}

\textit{Set of parameters}. Our selection rule consists in choosing
an estimator $\widehat{f}_{(h,\cP)}^{(n)}(x_{0})$ when the parameter
$(h,\cP)$ belongs at most to the set $\mH[ \mP]$ defined as follows.

Let $\mz>0$, $\tau(s)\in(0,1]$, $s=1,\ldots,d$, be fixed numbers and
let $\mathfrak{h}_I^{(I)}\in(0,1]^{\llvert I\rrvert }$, $I\in\cI_d$,
be fixed
multibandwidths. All these parameters will be chosen in accordance
with our procedure.

Set also $\lambda:=\sup_{I\in\cI_d} \{1\vee\lambda_{\llvert I\rrvert }^{(2q)} [\mathbf{K},\mz ] \}$ and $a:= \{
2\lambda
\sqrt{1+2q} \}^{-2}$, where constants $\lambda_s^{(q)}[\mathbf
{K},\mz
],   s\in\bN^*, q\geq1$, are given in Section~\ref{constants}. The
explicit expressions of $\lambda_s^{(q)}[\mathbf{K},\mz]$ are too
cumbersome and it is not convenient for us to present them right now.

For all $I\in\cI_d$ and all integer $m>0$ introduce
\begin{eqnarray*}
\mH_{m,1}^{(I)}&:=&\bigl\{h_I\in (0,1]^{|I|}\dvt v_m^{(I)}V_{\mh_I^{(I)}}\leq V_{h_I}
\leq v_{m-1}^{(I)}V_{\mh_I^{(I)}}\bigr\}\cap\prod_{i\in I}\biggl[\frac{1}{n},\bigl(v_m^{(I)}\bigr)^{-\mz}\mh_i^{(I)}\biggr],
\\[-2pt]
\mH_{m,2}^{(I)}&:=&\bigl\{h_I\in (0,1]^{|I|}\dvt v_m^{(I)}V_{\mathrm{max}}\leq V_{h_I}\leq v_{m-1}^{(I)}V_{\mathrm{max}}\bigr\}
\cap\prod_{i\in I}\biggl[\frac{1}{n},\bigl(v_m^{(I)}\bigr)^{-\mz}\mh_i^{(I)}\biggr],
\\[-2pt]
\mH^{(I)}&:=&\Biggl(\bigcup_{m=1}^{M_n(I)}\mH_{m,1}^{(I)}\Biggr)\bigcup\Biggl(\bigcup_{m=1}^{M_n(I)}\mH_{m,2}^{(I)}\Biggr),
\end{eqnarray*}
where $v_m^{(I)}:=2^{-m\tau(|I|)}$, $M_n(I)$ is the largest integer satisfying $v_{M_n(I)}^{(I)}[V_{\mh_I^{(I)}}\wedge V_{\mathrm{max}}]\geq
 \frac{\ln(n)}{an}$ and $M_n(I)\leq\log_2(n)$, and $V_{\mathrm{max}}$ is defined below.


Define finally
\[
\mH[ \mP]:= \bigl\{(h,\cP)\in(0,1]^d\times\mP\dvt h_I
\in\mH^{(I)}, \forall I\in\cP \bigr\}. %
\]

\textit{Extra parameters.} Let $\overline{\mH}$ and $\overline
{\mP}$
be arbitrary subsets of $(0,1]^d$ and $\mP$, respectively. The
selection rule (\ref{eq:selectionrule1})--(\ref{eq:selectionrule2})
below run over $\overline{\mH}[ \overline{\mP} ]:= (\overline
{\mH
}\times\overline{\mP} )\cap\mH[ \mP]$ and the reasons for
introducing these extra parameters are discussed in Remark~\ref
{selectionrule}. In particular, for measurability reasons, we will
always suppose that $\overline{\mH}$ is either a compact or a finite
subset of $(0,1]^d$.

Set $\Lambda_n(x_0):=3\lambda d^2 [ 2\overline
{G}_{n}(x_{0})
]^{d^2-1}$, where
\[
\overline{G}_{n}(x_0):=\sup_{(h,\cP)\in\overline{\mH}[ \overline
{\mP}
]}
\sup_{(\eta,\cP')\in\overline{\mH}[ \overline{\mP} ]}\sup_{I\in\cP
\circ\cP'} \bigl[2
\widetilde{G}_{h_I\vee\eta_I}(x_{0,I}) \bigr]. %
\]

Put also $V_{\mathrm{max}}:=\sup_{\cP\in\overline{\mP}}\inf_{I\in\cP
}V_{\mathfrak{h}
_I^{(I)}}$ and, for $(h,\cP)\in(0,1]^d\times\mP$,
\[
\delta(h,\cP):=\sup_{\cP'\in\overline{\mP}}\sup_{I\cap I'\in
\cP\circ\cP
'}
\biggl[\frac{V_{\mathfrak{h}_{I\cap I'}^{(I)}\vee\mathfrak
{h}_{I\cap
I'}^{(I')}}}{V_{h_{I\cap I'}}} \biggr]\vee \biggl[\frac{V_{\mathrm{max}}}{\inf_{I\in
\cP}V_{\mathfrak{h}_I^{(I)}}} \biggr]. %
\]

Define finally, for $(h,\cP)\in(0,1]^d\times\mP$,
\[
\widehat{\cU}_{(h,\cP)}(x_0):=\sqrt{\frac{ [ \overline{G}_{n}(x_0)
 ]^2 \{1\vee\ln\delta(h,\cP) \}}{nV(h,\cP)}},\qquad V(h,
\cP ):=\inf_{I\in\cP}V_{h_I}. %
\]

\subsection{Selection rule}
\label{sec:selectionrule}
For $(h,\cP)\in(0,1]^d\times\mP$ introduce
%
\begin{eqnarray}
\label{eq:selectionrule1} &&\widehat{\Delta}_{(h,\cP)}(x_0)
\nonumber
\\[-8pt]
\\[-8pt]
\nonumber
&&\quad:=\sup
_{(\eta,\cP')\in\overline
{\mH}[
\overline{\mP} ]} \bigl[\bigl\llvert \widehat{f}_{(h,\cP),(\eta,\cP
')}^{(n)}(x_{0})-
\widehat{f}_{(\eta,\cP')}^{(n)}(x_{0})\bigr\rrvert -\Lambda
_n(x_0) \bigl\{\widehat{\cU}_{(\eta,\cP')}(x_0)+
\widehat{\cU }_{(h,\cP
)}(x_0) \bigr\} \bigr]_+.\quad
\end{eqnarray}
Define finally $ (\widehat{h},\widehat{\cP} )$ satisfying
%
\begin{equation}
\label{eq:selectionrule2} \widehat{\Delta}_{(\widehat{h},\widehat{\cP})}(x_0)+2\Lambda
_n(x_0)\widehat{\cU}_{(\widehat{h},\widehat{\cP})}(x_0)=
\inf_{(h,\cP)\in
\overline{\mH}[ \overline{\mP} ]} \bigl[\widehat{\Delta}_{(h,\cP
)}(x_0)+2
\Lambda_n(x_0)\widehat{\cU}_{(h,\cP)}(x_0)
\bigr].
\end{equation}
The selected estimator is $\widehat{f}_{n}(x_{0}):=\widehat
{f}_{(\widehat{h},\widehat{\cP})}^{(n)}(x_{0})$.

Similarly to Section~2.1 in Lepski \cite{lepski:supnormlossdensityestimation} it is easy to show that $
(\widehat{h},\widehat{\cP} )$ is $X^{(n)}$-measurable and that
$
(\widehat{h},\widehat{\cP} )\in\overline{\mH}[ \overline{\mP
} ]$.
It follows that $\widehat{f}_{n}(x_{0})$ is also a $X^{(n)}$-measurable
random variable.

\begin{remark}
\label{selectionrule}
The necessity to introduce the extra parameters $\overline{\mH}$ and
$\overline{\mP}$ is dictated by
several reasons. The first one is computational namely the computation
of $\widehat{\Delta}_{(h,\cP)}(x_0)$ and $(\widehat{h},\widehat
{\cP})$. However,
the computational aspects of the choice of $\overline{\mP}$ and
$\overline{\mH}$ are quite different. Typically, $\overline{\mH}$ can
be chosen as an appropriate grid in $(0,1]^d$, for instance, dyadic one,
that is sufficient for proving adaptive properties of the proposed estimator.
The choice of $\overline{\mP}$ is much more delicate. The reason of
considering $\overline{\mP}$ instead of $\mP$ is explained by
the fact that the cardinality of $\mP$ grows exponentially with the
dimension $d$. Therefore, if $\overline{\mP}=\mP$, for large values of
$d$ our procedure is not practically feasible in view of huge amount of
comparisons to be done. In the latter case, the interest of our result
is theoretical. Note also that the best attainable trade-off between
approximation and stochastic errors depends heavily on both the number
of observations and the effective dimension $d(f)=\inf_{\cP\in\mP
(f)}\sup_{I\in\cP}\llvert I\rrvert $. Thus, if $d(f)$ is big the
corresponding independence structure does not bring a real improvement
of the estimation accuracy. So, in practice, $\overline{\mP}$ is chosen
to satisfy $\sup_{I\in\cP}\llvert I\rrvert \leq d_0$, $\forall\cP
\in
\overline{\mP} \setminus \{\overline{\varnothing} \}$. The
choice of the parameter $d_0$ (made by a statistician) is based on the
compromised between the sample size $n$, the desirable quality of
estimation and the number of computations. For instance, one can
consider $d_0=1$, that means that $\overline{\mP}$ contains two
elements, $\{\{1,\ldots,d\}\}$ and $\{\{1\},\ldots,\{d\}\}$. The latter
case corresponds to the observations having independent components and
it can be illustrated in Example~\ref{ex:densityestimation} below. On
the other hand, in the case of low dimension $d$, one can always take
$\overline{\mP}=\mP$, since if $d=2$, $\llvert \mP\rrvert =2$, $d=3$,
$\llvert \mP\rrvert =5$, $d=4$, $\llvert \mP\rrvert =12$, etc.

Other reasons are related to the possibility to consider various
problems arising in the framework of minimax and minimax adaptive
estimation and they will be discussed in detail in Sections~\ref{sec:minimaxresults} and \ref{sec:optimality}. Here, we only mention
that the choice $\overline{\mP}= \{\overline{\varnothing} \}$
allows to study the adaptive estimation of a multivariate density on
$\bR^d$ without taking into account eventual independence structure. We
would like to emphasize that the latter problem was not studied in the
literature.

At last the introduction of $\overline{\mP}$ allows to minimize the
assumptions imposed on the density to be estimated.
In particular, the oracle inequality corresponding to $\overline{\mP
}= \{\overline{\varnothing} \}$ is proved over the set of bounded
densities; see Corollary~\ref{cor:oracleinequality}.
\end{remark}

In spite of the fact that the construction of the proposed procedure
does not require any condition on the density $f$, the following
assumption will be used for computing its risk:
%
\begin{equation}
\label{eq:densityassumption} f\in\bF_d [\mathbf{f},\overline{\mP} ]:= \Bigl\{f\dvt
\sup_{\cP,\cP'\in\overline{\mP}}\sup_{I\in\cP\circ\cP'}\llVert
f_I\rrVert _{\infty
}\leq\mathbf{f}, \exists\cP\in\mP(f)\cap
\overline{\mP} \Bigr\},\qquad  0<\mathbf{f}<+\infty.
\end{equation}
Note that the considered class of densities is determined by $\overline
{\mP}$ and in particular, if $\overline{\varnothing}\in\overline
{\mP}$,
\begin{eqnarray*}
\bF_d [\mathbf{f},\mP ]&=& \Bigl\{f\dvt \sup_{I\in\cI
_d}
\llVert f_I\rrVert _{\infty}\leq\mathbf{f} \Bigr\}\subseteq
\bF_d [\mathbf {f},\overline{\mP} ],
\\
\bF_d \bigl[\mathbf{f}, \{\overline{\varnothing} \} \bigr]&=& \bigl\{
f\dvt \llVert f\rrVert _{\infty}\leq\mathbf{f} \bigr\}, \qquad \bF_d
\bigl[\mathbf{f}, \{\cP \} \bigr]= \Bigl\{f\dvt \sup_{I\in\cP
}
\llVert f_I\rrVert _{\infty}\leq\mathbf{f}, \cP\in\mP(f)
\Bigr\}.
\end{eqnarray*}

\subsection{Oracle-type inequality}
\label{sec:oracleinequality}
For $I\in\cI_d$ and $(h,\eta)\in(0,1]^d\times[0,1]^d$ introduce
\[
\cB_{h_I,\eta_I}(x_{0,I}):=\int_{\bR^{\llvert I\rrvert }}K^{(I)}(u)
\bigl[f_I \bigl(x_{0,I}+(h_I\vee
\eta_I)u \bigr)-f_I (x_{0,I}+
\eta_I u ) \bigr]\,\rd u, %
\]
where here and later $y_Ix_I$ denotes the coordinate-vise product of
$y_I,x_I\in\bR^{\llvert I\rrvert }$.

For $(h,\cP)\in(0,1]^d\times\mP$ define $\cB_{(h,\cP)}(x_0):=\sup_{\cP
'\in\overline{\mP}}\sup_{I\in\cP\circ\cP'}\sup_{\eta\in
[0,1]^d}\llvert \cB
_{h_I,\eta_I}(x_{0,I})\rrvert $.

Introduce finally, if exists $\cP\in\mP(f)\cap\overline{\mP}$,
\[
\mR_n(f):=\inf_{(h,\cP)\in\overline{\mH}[ \overline{\mP} ]\dvt\cP
\in\mP
(f)} \biggl[
\cB_{(h,\cP)}(x_0)+\sqrt{\frac{1\vee\ln\delta(h,\cP
)}{nV(h,\cP)}} \biggr].
\]
The quantity $\mR_n(f)$ can be viewed as the optimal trade-off between
approximation and stochastic errors provided by estimators involved in
the selection rule.

\begin{theorem}
\label{theo:oracleinequality}
Let $\overline{\mH}\subseteq(0,1]^d$ and $\overline{\mP}\subseteq
\mP$
be arbitrary subsets such that $\overline{\mH}[ \overline{\mP} ]$ is
non-empty.

Then for any $0<\mathbf{f}<+\infty$, any $q\geq1$ and any integer
$n\geq3$:
%
\begin{equation}
\label{eq:oracleinequality} \cR_n^{(q)} [\widehat{f}_n,f ]
\leq\alpha_1\mR _n(f)+\alpha_2
[nV_{\mathrm{max}} ]^{-{1}/{2}}\qquad \forall f\in\bF_d [\mathbf {f},
\overline{\mP} ],
\end{equation}
where $\alpha_1:=\alpha_1(q,d,\mathbf{K},\mathbf{f})$ and $\alpha
_2:=\alpha_2(q,d,\mathbf{K},\mathbf{f})$ are given in the proof of
the theorem.
\end{theorem}

Considering the case $\overline{\mP}= \{\overline{\varnothing
}
\}$ and noting $\mH=\mH^{(\overline{\varnothing})}$ we come to the
following consequence of Theorem~\ref{theo:oracleinequality}.

\begin{corollary}
\label{cor:oracleinequality}
Let assumptions of Theorem~\ref{theo:oracleinequality} be fulfilled.
Then, for all densities $f$ such that $\llVert f\rrVert _{\infty}\leq
\mathbf{f}$,
%
\begin{equation}
\label{eq:oracleinequalitywithoutindependence} \cR_n^{(q)} [\widehat{f}_n,f ]
\leq\alpha_1\inf_{h\in
\overline
{\mH}\cap\mH} \biggl[ \sup
_{\eta\in[0,1]^d}\bigl\llvert \cB_{h,\eta
}(x_{0})\bigr
\rrvert +\sqrt{\frac{1\vee\ln ({V_{\mathfrak
{h}}}/{V_{h}}
)}{nV_{h}}} \biggr]+\alpha_2
[nV_{\mathfrak{h}} ]^{-{1}/{2}}.
\end{equation}
\end{corollary}

Looking at the assertion of Theorem~\ref{theo:oracleinequality} and its
Corollary~\ref{cor:oracleinequality} it is not clear what can be gained
by taking into account eventual independence structure. This issue will
be scrutinized in Section~\ref{minimax-adaptive-estimation}, but some
conclusions can be deduced directly from the latter results. Consider
the following example.

\begin{example}
\label{ex:densityestimation}
For any $t\in\bR$, put
\[
f(t)=\tfrac{64}{15} \bigl\{4t\mathbf{1}_{[0,1/8)}(t)+ \bigl(
\tfrac
{3}{4}-2t \bigr)\mathbf{1}_{(1/8,1/4]}(t)+\tfrac{1}{4}
\mathbf {1}_{(1/4,3/4]}(t)+ (1-t )\mathbf{1}_{(3/4,1]}(t) \bigr\},
\]
and define $f_d(x)=\prod_{i=1}^d f(x_i)$, $x\in\bR^d$. It is easily
seen that $f_d$ is a probability density and the goal is to estimate
$f(x_0)$, $x_0\in(3/8,7/8)^d$.

Choose $\mathfrak{h}=(1,\ldots,1)$, $h=(1/4,\ldots,1/4)$ and let
$\overline{\mH
}=\{\mathfrak{h},h\}$. Put $\cP_1=\{\{1,\ldots,d\}\}$, $\cP_2=\{\{
1\},\ldots,\{
d\}\}$ and let $\overline{\mP}=\{\cP_1,\cP_2\}$. Since, in this case,
$\overline{\mH}\times\overline{\mP}$ contains $4$ elements, our
estimator can be computed in a reasonable time.

Moreover, in accordance with the oracle-type inequality proved in
Theorem~\ref{theo:oracleinequality}, the accuracy provided by the selected estimator is
proportional to $\sqrt{[4\ln (4 )]/n}$. On the other hand, the
pointwise risk of the kernel estimator with optimally chosen bandwidth
and kernel is proportional to $\sqrt{[d4^d\ln (4 )]/n}$ if the
independence structure is not taken into account. As we see, the
adaptation to eventual independence structure can lead to significant
improvement of the constant. This shows that the proposed methodology
has an interest beyond derivation of minimax rates, which is the
subject of the next section.
\end{example}

\section{Minimax and adaptive minimax pointwise estimation}
\label{minimax-adaptive-estimation}

In this section, we provide with minimax and adaptive minimax
estimation over a scale of anisotropic Nikolskii classes.

\subsection{Anisotropic Nikolskii densities classes with independence
structure}
\label{sec:nikolskiiclasses}

Let $ \{e_1,\ldots,e_s \}$ denote the canonical basis in
$\bR
^s, s\in\bN^*$.

\begin{definition}
Let $p= (p_1,\ldots,p_s), p_i\in[1,\infty]$, $\beta= (\beta
_1,\ldots,\beta_s), \beta_i>0$ and $L= (L_1,\break\ldots,L_s)$, $L_i>0$. A function
$f\dvtx\bR^s\rightarrow\bR$ belongs to the anisotropic Nikolskii class
$\bN
_{p,s}(\beta,L)$~if
\begin{eqnarray*}
\mathrm{(i)}&&\quad \bigl\llVert D_i^{k} f\bigr\rrVert
_{p_i}\leq L_i\qquad \forall k=\overline {0, \lfloor
\beta_i \rfloor}, \forall i=\overline{1,s};
\\
\mathrm{(ii)}&&\quad \bigl\llVert D_i^{ \lfloor\beta_i \rfloor} f(\cdot
+te_i)-D_i^{ \lfloor\beta_i \rfloor} f(\cdot)\bigr\rrVert
_{p_i}\leq L_i\llvert t\rrvert ^{\beta_i- \lfloor\beta_i
\rfloor
} \qquad\forall t
\in\bR, \forall i=\overline{1,s}.
\end{eqnarray*}
Here, $D_i^{k} f$ denotes the kth order partial derivate of $f$ with
respect to the variable $t_i$, and $ \lfloor\beta_i
\rfloor$
is the largest integer strictly less than $\beta_i$.
\end{definition}

The following\vspace*{1pt} collection $ \{N_{p,d}^* (\beta,L,\cP
)
\}_{\cP}$ was introduced in Lepski \cite{lepski:supnormlossdensityestimation} in order to take into account the
smoothness of the underlying density and its eventual independence
structure simultaneously.
\[
\label{eq:adaptiveminimaxclass1} N_{p,d}^* (\beta,L,\cP ):= \biggl\{f\in
\bN_{p,d}^* (\beta,L )\dvt f\geq0, \int{f}=1, f(x)=\prod
_{I\in\cP}f_I(x_I), \forall x\in
\bR^d \biggr\},
\]
where $f\in\bN_{p,d}^* (\beta,L )$ means that
%
\begin{equation}
\label{eq:marginalenikolskii1} f_I\in\bN_{p_I,\llvert I\rrvert }
(\beta_I,L_I)\qquad
\forall I\in\cI_d.
\end{equation}

We remark that this collection of functional classes was used in the
case of adaptive estimation, that is, when the partition $\cP\in\mP$ is
unknown. However, when the minimax estimation is considered ($\cP$ is
fixed), we do not need that condition (\ref{eq:marginalenikolskii1})
holds for any $I\in\cI_d$. It suffices to consider only $I$ belonging
to $\cP$, and we come to the following definition.

\begin{definition}[(Minimax estimation)]
\label{def:nikolskiidensityminimax}
Let $p= (p_1,\ldots,p_d), p_i\in[1,\infty]$, $\beta= (\beta
_1,\ldots,\beta_d),\break \beta_i>0$, $L= (L_1,\ldots,L_d)$, $L_i>0$ and $\cP\in
\mP$.
A probability density $f\dvtx\bR^d\rightarrow\bR_+$ belongs to the class
$N_{p,d} (\beta,L,\cP )$ if
%
\begin{eqnarray}
\label{eq:marginalenikolskii} f(x)=\prod_{I\in\cP}f_I(x_I)\qquad
\forall x\in\bR^d, \qquad f_I\in\bN _{p_I,\llvert I\rrvert }(
\beta_I,L_I) \qquad\forall I\in\cP.
\end{eqnarray}
\end{definition}

Let us now come back to the adaptive estimation. As it was discussed in
Remark~\ref{selectionrule}, the adaptation is not necessarily
considered with respect to $\mP$. If $\overline{\mP}\subset\mP$ is used
instead of $\mP$, the assumption (\ref{eq:marginalenikolskii1}) is too
restrictive and can be weakened in the following way.

Denote $\overline{\mP}^*:= \{\cP\circ\cP^{\prime}\dvt \cP, \cP
^{\prime}\in\overline{\mP} \}$ and $\overline{\cI
}_d^*:= \{I\in
\cI_d\dvt \exists\cP\in\overline{\mP}^*, I\in\cP \}$.

\begin{definition}[(Adaptive estimation)]
\label{def:nikolskiidensityadaptive}
Let $\overline{\mP}\subseteq\mP$ and $ (\beta,p,\cP
)\in
(0,+\infty )^d\times[1,\infty]^d\times\overline{\mP}$ be
fixed. A
probability density $f\dvtx\bR^d\rightarrow\bR_+$ belongs to the class
$\overline{N}_{p,d} (\beta,L,\cP )$ if
%
\begin{equation}
\label{eq:marginalenikolskii2} f(x)=\prod_{I\in\cP}f_I(x_I)\qquad
\forall x\in\bR^d;\qquad  f_I\in\bN _{p_I,\llvert I\rrvert }(
\beta_I,L_I)\qquad \forall I\in\overline{
\cI}_d^*.
\end{equation}
\end{definition}

Some remarks are in order.

(1) We note that if $\overline{\mP}=\mP$, then $\overline
{N}_{p,d}
(\beta,L,\cP )=N_{p,d}^* (\beta,L,\cP )$, but for some
$\overline{\mP}\subset\mP$, one has $N_{p,d}^* (\beta,L,\cP
 )$
$\subset\overline{N}_{p,d} (\beta,L,\cP )$. The latter
inclusion shows that the condition (\ref{eq:marginalenikolskii2}) is
weaker than $f\in N_{p,d}^* (\beta,L,\cP )$. In
particular, if
$\overline{\mP}= \{\overline{\varnothing} \}$, then
$\overline
{N}_{p,d} (\beta,L,\overline{\varnothing} )= \{f\in\bN
_{p,d}(\beta,L)\dvt f\geq0, \int f=1 \}\supset N_{p,d}^* (\beta,L,\overline{\varnothing} )$.

(2) Note that if $ \overline{\mP}= \{\cP \}$, then
$\overline
{N}_{p,d} (\beta,L,\cP )$ coincides with the class
$N_{p,d} (\beta,L,\cP )$ used for minimax estimation. But
$\overline{N}_{p,d} (\beta,L,\cP )\subset N_{p,d}
(\beta,L,\cP )$ for all $\cP\in\overline{\mP}$ for any other
choices of~$ \overline{\mP}$.

\subsection{Minimax results}
\label{sec:minimaxresults}
For $ (\beta,p,\cP )\in (0,+\infty )^d\times
[1,\infty
]^d\times\mP$ define
%
\begin{eqnarray}
r:=r (\beta,p,\cP )&=&\inf_{I\in\cP}\gamma_I (\beta,p
), \qquad \gamma_I:=\gamma_I (\beta,p )=\frac{1-\sum_{i\in
I}{1}/{(\beta_i p_i)}}{\sum_{i\in I}{1}/{\beta_i}},\qquad
I\in\cP; \nonumber
\\
\label{eq:minimaxrate} \varphi_n (\beta,p,\cP )&:= &\biggl(\frac{1}{n}
\biggr)^{{r}/{(2r+1)}}, \qquad \rho_n (\beta,p,\cP ):=\mathbf{1}_{
\{r\leq0 \}}+
\varphi_n (\beta,p,\cP )\mathbf {1}_{ \{
r> 0 \}}.
\end{eqnarray}
As it will follow from Theorems \ref{theo:minimaxlowerbound1} and \ref
{theo:minimaxupperbound} below $\varphi_n (\beta,p,\cP )$ is
the minimax rate of convergence on $N_{p,d}(\beta,L,\cP)$. Hence,
similarly to the standard representation of minimax rates, the
parameter~$r$ can be interpreted as a smoothness index corresponding to
the independence structure.

\begin{theorem}
\label{theo:minimaxlowerbound1}
$\forall (\beta,p,\cP )\in (0,+\infty
)^d\times
[1,\infty]^d\times\mP$, $\forall L\in(0,\infty)^d$, $\exists c>0$:
\[
\liminf_{n\rightarrow+\infty} \Bigl\{\rho_n^{-1} (
\beta,p,\cP )\inf_{\widetilde{f}_n}\cR_n^{(q)} \bigl[
\widetilde {f}_n,N_{p,d}(\beta,L,\cP) \bigr] \Bigr\}\geq c,
\]
where infimum is taken over all possible estimators.
\end{theorem}

Note that the assertion of Theorem~\ref{theo:minimaxlowerbound1} will
be deduced from more general result established in Proposition~\ref
{theo:minimaxlowerbound} below. It is also important to emphasize that
if $r\leq0$ there is no uniformly consistent estimator for the
considered problem and, to the best of our knowledge, this fact was not
known before. Let us provide an example with a density for which $r<0$.

\begin{example} Suppose that $d=1$ and, therefore, $\cP=\overline
{\varnothing}$ (no independence structure). For any $x\in\bR$, put
\[
g(x)=\mathbf{1}_{\{0\}}(x)+\frac{1}{2\sqrt{x}}\mathbf{1}_{(0,1]}(x).
\]
Some straightforward computations allows us to assert that $g\notin
N_{p,1} (\beta,L,\overline{\varnothing} )$, $\forall L>0$, if
$p\beta\geq1$ (i.e., $r\geq0$), and that $g\in N_{1,1}
(1/2,L,\overline{\varnothing} )$ for some $L>0$ ($p=1$, $\beta
=1/2$). Thus, in this case, one has $r<0$.
\end{example}

Our goal now is to show that $\varphi_n (\beta,p,\cP )$
is the
minimax rate of convergence on $N_{p,d}(\beta,L,\cP)$ and that a
minimax estimator belongs to the collection $\mF[ \mP]$. In fact, we
prove that the minimax estimator is $\widehat{f}_{(\mathbf{h},\cP
)}^{(n)}$ with properly chosen kernel $\mathbf{K}$ and bandwidth~$\mathbf{h}$.

For a given integer $l\geq2$ and a given symmetric Lipschitz function
$u\dvtx\bR\rightarrow\bR$ satisfying $\operatorname{supp}(u)\subseteq[-1/(2l),1/(2l)]$ and
$\int_{\bR} u(y)\,\mathrm{d}y=1$ set
%
\begin{equation}
\label{eq:orthogonalkernel} u_l(z):=\sum_{i=1}^{l}
\pmatrix{ l
\cr
i } %
(-1)^{i+1}\frac{1}{i}u \biggl(
\frac{z}{i} \biggr), \qquad z\in\bR.
\end{equation}
Furthermore, we use $\mathbf{K}\equiv u_l$ in the definition of
estimators collection $\mF[ \mP]$.

The relation of kernel $u_l$ to anisotropic Nikolskii classes is
discussed in Kerkyacharian, Lepski and Picard \cite{kerk1}. In
particular, it was shown that
%
\begin{equation}
\label{eq:kernelorthogonality} \int_{\bR}\mathbf{K}(z)\,\rd z=1,\qquad \int
_{\bR}z^k\mathbf{K}(z)\,\rd z=0\qquad \forall k=1,
\ldots,l-1.
\end{equation}

Choose finally $\mathbf{h}=(\mathbf{h}_1,\ldots,\mathbf{h}_d)$, where
\[
\mathbf{h}_i=n^{-({\gamma_I(\beta,p)}/{(2\gamma_I(\beta,p)+1)})
({1}/{\beta_i(I)})},\qquad i\in I, I\in\cP. %
\]
Here,
\[
\beta_i(I):=\kappa(I)\beta_i\kappa_i^{-1}(I),\qquad
\kappa(I):=1-\sum_{k\in I} (\beta_kp_k
)^{-1},\qquad \kappa_i(I):=1-\sum_{k\in
I}
\bigl(p_k^{-1}-p_i^{-1} \bigr)
\beta_k^{-1}. %
\]

\begin{theorem}
\label{theo:minimaxupperbound}
For all $ (\beta,p,\cP )\in(0,l]^d\times[1,\infty
]^d\times\mP$
such that $r (\beta,p,\cP )>0$ and all $L\in(0,\infty)^d$
\[
\limsup_{n\rightarrow+\infty} \bigl\{\varphi_n^{-1} (
\beta,p,\cP )\cR_n^{(q)} \bigl[\widehat{f}_{(\mathbf{h},\cP
)}^{(n)},N_{p,d}(
\beta,L,\cP) \bigr] \bigr\}<\infty. %
\]
\end{theorem}

To get the statement of this theorem, we apply Theorem~\ref
{theo:oracleinequality} with $\overline{\mP}=\{\cP\}$ and $\overline
{\mH
}=\{\mathbf{h}\}$. In view of the embedding theorem for anisotropic
Nikolskii classes (formulated in the proof of Lemma~\ref
{lem:biasupperbound} and available when $r (\beta,p,\cP )>0$),
there exists a number $\mathbf{f}:=\mathbf{f}(\beta,p)>0$ such that
$N_{p,d}(\beta,L,\cP)\subseteq\bF_d [\mathbf{f}, \{\cP
 \}
 ]$. It makes possible the application of Theorem~\ref
{theo:oracleinequality}.

Let us briefly discuss several consequences of Theorems \ref
{theo:minimaxlowerbound1} and \ref{theo:minimaxupperbound}.
First, if $\cP=\overline{\varnothing}$, we obtain the minimax rate on
the anisotropic Nikolskii class $\bN_{p,d}(\beta,L)$. In particular, if
$p_i=+\infty$, $i=\overline{1,d}$, we find the minimax rate on the
anisotropic H\"older class $\bH_d(\beta,L)$ given in (\ref
{eq:multivariateminimaxrates}). If $d=1$, then our results coincide
with those presented in (\ref{eq:univariateminimaxrates}).

Next, in view of Theorem~\ref{theo:minimaxlowerbound1} there is no
consistent estimator for $f(x_0)$ on $\bN_{p,d}(\beta,L)$ if $r
(\beta,p,\overline{\varnothing} )\leq0$. On the other hand, if
$f\in N_{p,d}(\beta,L,\cP)$ and $r(\beta,p,\cP)>0$, then such estimator
for $f(x_0)$ does exist in view of Theorem~\ref{theo:minimaxupperbound}
even if $r (\beta,p,\overline{\varnothing} )<0$.

Note also that the condition $r (\beta,p,\overline{\varnothing
}
)>0$ is sufficient to find a consistent estimator on each functional
class $N_{p,d}(\beta,L,\cP)$, $\cP\in\mP$, and that the same condition
is necessary for the estimation over $N_{p,d}(\beta,L,\varnothing)$. It
allows us to compare the influence of the independence structure on the
accuracy of estimation. For example, we see that
\[
\varphi_n \bigl(\bH_d(\beta,L) \bigr)\gg
\varphi_n (\beta,p,\cP ), \qquad p_i=\infty, i=\overline{1,d}.
\]

We conclude that the existence of an independence structure improves
significantly the accuracy of estimation.

We finish this section with the result being a refinement of Theorem~\ref{theo:minimaxlowerbound1}.

\begin{proposition}
\label{theo:minimaxlowerbound}
$\forall (\beta,p,\cP )\in (0,+\infty
)^d\times
[1,\infty]^d\times\mP$, $\forall L\in(0,\infty)^d$, $\exists c>0$:
\begin{eqnarray*}
\liminf_{n\rightarrow+\infty} \Bigl\{\rho_n^{-1} (
\beta,p,\cP )\inf_{\widetilde{f}_n}\cR_n^{(q)} \bigl[
\widetilde {f}_n,N_{p,d}^*(\beta,L,\cP) \bigr] \Bigr\}\geq
c,
\end{eqnarray*}
where infimum is taken over all possible estimators.
\end{proposition}

\begin{remark}
\label{rem:minimaxrate}
Recall (see Section~\ref{sec:nikolskiiclasses}) that $N_{p,d}^*(\beta,L,\cP)\subseteq\overline{N}_{p,d} (\beta,L,\cP
)\subseteq
N_{p,d} (\beta,L,\cP )$. Hence, the statement of Theorem~\ref
{theo:minimaxupperbound} remains true if one replaces $N_{p,d}
(\beta,L,\cP )$ by $\overline{N}_{p,d} (\beta,L,\cP
 )$,
$\overline{\mP}\subseteq\mP$. Thus, Proposition~\ref
{theo:minimaxlowerbound} together with Theorem~\ref
{theo:minimaxupperbound} allows us to assert that $\rho_n (\beta,p,\cP )$ is the minimax rate of convergence on $\overline
{N}_{p,d} (\beta,L,\cP )$.
\end{remark}

\subsection{Adaptive estimation}
\label{sec:adaptiveminimaxresults}

\subsubsection{Adaptive estimation. Upper bound}

Let $\overline{\mP}\subseteq\mP$, such that $\overline{\varnothing
}\in
\overline{\mP}$, be fixed. Denote $d(\cP):=\sup_{I\in\cP}|I|$,
$\cP\in
\overline{\mP}$, and $\overline{d}:=\inf_{\cP\in\overline{\mP
}}\,\mathrm{d}(\cP)$.

Set $\beta_i^{(\mathrm{max})}=\beta_{\mathrm{max}}>(d-\overline{d})/2$,
$p_i^{(\mathrm{max})}=+\infty, i=\overline{1,d}$, and suppose additionally
that $l\geq2\vee\beta_{\mathrm{max}}$.
Choose $\mathbf{K}\equiv u_l$, $\mz:=\frac{1}{2\beta_{\mathrm{max}}}$ and
$\tau
(s)$, $s=1,\ldots,d$, satisfying
\[
\tau(s):=2\beta_{\mathrm{max}}/(2\beta_{\mathrm{max}}+\overline{d}).
\]

Let $\overline{\mH}$ be the dyadic grid in $(0,1]^d$ and let
$\mathfrak{h}
_I^{(I)}$, $I\in\cI_d$, be the projection on the dyadic grid in
$(0,1]^{|I|}$ of the multibandwith $\rh_I^{(I)}$ given by
%
\begin{equation}
\label{eq:maximalebandwidth}
\rh_i^{(I)}:=n^{-1/(2\beta_{\mathrm{max}}+\overline{d})},\qquad i\in I.
\end{equation}

Consider the estimator $\widehat{f}_n(x_0)$ defined by the selection
rule (\ref{eq:selectionrule1})--(\ref{eq:selectionrule2}), in
Section~\ref{sec:selectionrule}.

For $ (\beta,p,\cP )\in (0,\beta_{\mathrm{max}}
]^d\times
[1,\infty]^d\times\overline{\mP}$ introduce
%
\begin{equation}
\label{eq:adaptiverate} \psi_n (\beta,p,\cP ):=\cases{ %
\displaystyle \biggl(
\frac{\ln(n)}{n} \biggr)^{{r}/{(2r+1)}}, &\quad
$r:=r (\beta,p,\cP )<r_{\mathrm{max}},$
\vspace*{2pt}
\cr
\displaystyle\biggl(\frac{1}{n} \biggr)^{{r_{\mathrm{max}}}/{(2r_{\mathrm{max}}+1)}},
 &\quad $r:=r (
\beta,p,\cP )=r_{\mathrm{max}},$}\qquad  r_{\mathrm{max}}:=\frac{\beta
_{\mathrm{max}}}{\overline{d}}.
\end{equation}

\begin{theorem}
\label{theo:adaptiveupperbound}
For any $ (\beta,p )\in (0,\beta_{\mathrm{max}}
]^d\times
[1,\infty]^d$ such that $r (\beta,p,\overline{\varnothing} )>0$,
any $\cP\in\overline{\mP}$ and any $L\in(0,\infty)^d$
\[
\limsup_{n\rightarrow+\infty} \bigl\{\psi_n^{-1} (
\beta,p,\cP )\cR_n^{(q)} \bigl[\widehat{f}_n,
\overline{N}_{p,d}(\beta,L,\cP ) \bigr] \bigr\}<\infty. %
\]
\end{theorem}

Similarly to Theorem~\ref{theo:minimaxupperbound}, the proof of Theorem~\ref{theo:adaptiveupperbound} is mostly based on the result of Theorem~\ref{theo:oracleinequality}. The application of Theorem~\ref
{theo:oracleinequality} is possible because $\overline{N}_{p,d}(\beta,L,\cP)\subseteq\bF_d [\mathbf{f},\overline{\mP}  ]$
for some
$\mathbf{f}:=\mathbf{f}(\beta,p)>0$ that is guaranteed by the condition
$r (\beta,p,\overline{\varnothing} )>0$.\vspace*{1pt}

We would like to emphasize that the construction of $\widehat
{f}_n(x_0)$ does not involved the knowledge of the parameters $(\beta,L,p,\cP)$. Using the modern statistical language, one can say that
$\widehat{f}_n(x_0)$ is fully adaptative.

Note, however, that the precision $\psi_n(\beta,p,\cP)$ given by this
estimator does not coincide with minimax rate of convergence $\varphi
_n(\beta,p,\cP)$ whenever $r\neq r_{\mathrm{max}}$. In the next section, we
prove that $\psi_n(\beta,p,\cP)$ found in Theorem~\ref
{theo:adaptiveupperbound} is an optimal payment for adaptation.

\subsubsection{Adaptive estimation. Criterion of optimality}
\label{sec:optimality}
Let $ \{\Sigma_{(\alpha,b)}, (\alpha,b)\in\cA\times\mB
 \}$ be
the scale of functional classes where $\cA\subset\bR^m$ is a
$(m)$-dimensional manifold and $\mB$ is a finite set. Recall that the
family $\varPsi= \{\varPsi_n(\alpha,b), (\alpha,b)\in\cA
\times\mB
 \}$ of normalizations is called admissible if there exists an
estimator $\widehat{f}_{\Psi}$ such that
%
\begin{equation}
\label{eq:adaptiv1} \limsup_{n\rightarrow+\infty} \bigl\{\varPsi_n^{-1}(
\alpha,b)\cR _n^{(q)} [\widehat{f}_{\Psi},
\Sigma_{(\alpha,b)} ] \bigr\} <+\infty\qquad \forall(\alpha,b)\in\cA\times\mB.
\end{equation}
The estimator $\widehat{f}_{\varPsi}$ is called $\varPsi$-adaptive.

In the considered problem, $\alpha=(\beta,p)$, $b=\cP$ and
\[
\cA= \bigl\{ (\beta,p )\in (0,\beta_{\mathrm{max}} ]^d\times [1,
\infty]^d\dvt r (\beta,p,\overline{\varnothing} )>0 \bigr\},\qquad \mB=
\overline{\mP}. %
\]
As it follows from Theorem~\ref{theo:adaptiveupperbound} $\psi
_n(\beta,p,\cP)$ is an admissible family of normalizations and the estimator
$\widehat{f}_n$ is $\psi_n$-adaptive.

Let $\varPsi= \{\varPsi_n(\alpha,b)>0, (\alpha,b)\in\cA
\times\mB
 \}$ and $\widetilde{\varPsi}= \{\widetilde{\varPsi
}_n(\alpha,b)>0, (\alpha,b)\in\cA\times\mB \}$ be arbitrary families of
normalizations and put
\[
\Upsilon_n(\alpha,b):=\frac{\widetilde{\varPsi}_n (\alpha,b
)}{\varPsi_n (\alpha,b )},\qquad \Upsilon_n(
\alpha):=\inf_{b\in
\mB}\Upsilon_n(\alpha,b).
\]
Define the set $\cA^{(0)} [\widetilde{\varPsi}/\varPsi
]\subseteq
\cA$ as follows:
\[
\cA^{(0)} [\widetilde{\varPsi}/\varPsi ]:= \Bigl\{\alpha \in\cA\dvt
\lim_{n\to\infty}\Upsilon_n(\alpha)=0 \Bigr\}. %
\]
The set $\cA^{(0)} [\widetilde{\varPsi}/\varPsi ]$ can be viewed
as the set where the family $\widetilde{\varPsi}$ ``outperforms'' the
family $\varPsi$.
For any $b\in\mB$, introduce
\[
\cA_b^{(\infty)} [\widetilde{\varPsi}/\varPsi ]:= \Bigl\{
\alpha\in \cA\dvt \lim_{n\to\infty}\Upsilon_n(
\alpha_0)\Upsilon_n(\alpha,b)=\infty, \forall
\alpha_0\in\cA^{(0)} [\widetilde{\varPsi}/\varPsi ] \Bigr
\}. %
\]
Remark first that the set $\cA_b^{(\infty)} [\widetilde{\varPsi
}/\varPsi ]$ is the set where the family $\varPsi$ ``outperforms'' the
family $\widetilde{\varPsi}$. Moreover, the ``gain'' provided by
$\varPsi
$ with respect to $\widetilde{\varPsi}$ on $\cA_b^{(\infty)}
[\widetilde{\varPsi}/\varPsi ]$ is much larger than its ``loss'' on
$\cA^{(0)} [\widetilde{\varPsi}/\varPsi ]$.

The idea led to the criterion of optimality formulated below is to say
that $\varPsi$ is ``better'' than $\widetilde{\varPsi}$ if there exists
$b\in\mB$ for which the set $\cA_b^{(\infty)} [\widetilde
{\varPsi
}/\varPsi ]$ is much more ``massive'' than $\cA^{(0)}
[\widetilde
{\varPsi}/\varPsi ]$.

\begin{definition}
\label{def:adaptiveoptimality}
\textup{(I)} A family of normalizations $\varPsi$ is called adaptive rate of
convergence if
\begin{longlist}[1.]
\item[1.]$\varPsi$ is an admissible family of normalizations;
\item[2.] for any admissible family of normalizations $\widetilde{\varPsi}$
satisfying $\cA^{(0)} [\widetilde{\varPsi}/\varPsi ]\neq
\varnothing$
\begin{itemize}
\item$\cA^{(0)} [\widetilde{\varPsi}/\varPsi ]$ is
contained in
a $(m-1)$-dimensional manifold,
\item there exists $b\in\mB$ such that $\cA_b^{(\infty)}
[\widetilde
{\varPsi}/\varPsi ]$ contains an open set of $\cA$.
\end{itemize}
\end{longlist}
\textup{(II)} If $\varPsi$ is an adaptive rate of convergence, then $\widehat
{f}_{\Psi}$ satisfying (\ref{eq:adaptiv1}) is called rate adaptive estimator.
\end{definition}
The aforementioned definition is inspired by Klutchnikoff's criterion;
see Klutchnikoff \cite{klutch}. Indeed if $\operatorname{card}(\mB)=1$ the both
definitions coincide.

\begin{theorem}
\label{theo:adaptiveoptimality} \textup{(i)} We can find no optimal rate
adaptive estimator (satisfying (\ref{eq:adaptivoptimal}) in
Section~\ref{sec:introduction}) over the scale
\[
\bigl\{ \overline{N}_{p,d}(\beta,L,\cP), (\beta,p,L,\cP)\in\mA \bigr
\}, %
\]
whenever $\mA\subseteq \{ (\beta,p,L,\cP )\in
(0,\beta
_{\mathrm{max}} ]^d\times[1,\infty]^d\times(0,\infty)^d\times\overline
{\mP}\dvt
r (\beta,p,\overline{\varnothing} )>0 \}$ contains at least
two elements $ (\beta,p,L,\cP )$ and $ (\beta',p',L',\cP
'
)$ such that $r(\beta,p,\cP)\neq r(\beta',p',\cP')$.

\textup{(ii)} $\widehat{f}_n(x_0)$ is rate adaptive estimator of $ f(x_0)$ and
$\psi_n$ is the adaptive rate of convergence, in the sense of
Definition~\ref{def:adaptiveoptimality}, over the scale
\[
\bigl\{ \overline{N}_{p,d}(\beta,L,\cP), (\beta,p,L,\cP )\in (0,
\beta_{\mathrm{max}} ]^d\times[1,\infty]^d\times(0,\infty
)^d\times \overline{\mP}, r (\beta,p,\overline{\varnothing} )>0 \bigr
\}. %
\]
\end{theorem}

It is important to emphasize that our results cover a large class of
problems in the framework of pointwise density estimation.

In particular, if $\overline{\mP}= \{\overline{\varnothing}
\}$,
we deduce that $\widehat{f}_n(x_0)$ is rate adaptive estimator of
$f(x_0)$ over
\[
\bigl\{\overline{N}_{p,d}(\beta,L,\overline{\varnothing}), (\beta,p,L)
\in (0,\beta_{\mathrm{max}} ]^d\times[1,\infty]^d\times
(0,\infty )^d, r (\beta,p,\overline{\varnothing} )>0 \bigr\}.
\]
The adaptive rate of convergence for this problem is given by
\[
\psi_n (\beta,p,\overline{\varnothing} ):=\cases{ %
 \displaystyle\biggl(\frac{\ln(n)}{n} \biggr)^{{r}/{(2r+1)}}, &\quad
 $\displaystyle (\beta,p )\neq \bigl(
\beta^{(\mathrm{max})},p^{(\mathrm{max})} \bigr),
r:=\frac{1-\sum_{i=1}^d{1}/{(\beta_i p_i)}}
{\sum_{i=1}^d{1}/{\beta_i}},$ \vspace*{2pt}
\cr
\displaystyle\biggl(\frac{1}{n} \biggr)^{{r_{\mathrm{max}}}/
{(2r_{\mathrm{max}}+1)}}, &\quad $\displaystyle (\beta,p )= \bigl(
\beta^{(\mathrm{max})},p^{(\mathrm{max})} \bigr), r_{\mathrm{max}}:=
\frac
{\beta_{\mathrm{max}}}{d}.$} %
\]
To the best of our knowledge, the latter result is new. It is precise
and generalizes the results of Butucea \cite{butucea:pointwisedensityadapt} ($d=1$) and Comte and Lacour \cite{comte:deconvolution} for the deconvolution model when the noise
variable is equal to zero.

Another interesting fact is related to the set of ``nuisance'' parameters
where the adaptive rate of convergence $\psi_n (\beta,p,\cP
 )$
coincides with the minimax one. In all known for us problems of
pointwise adaptive estimation this set contains a single element.
However, as it follows from Theorem~\ref{theo:adaptiveoptimality}, this
set may contain several elements. Indeed, if, for instance, $d=4$ and
$\overline{\mP}= \{\cP_1,\cP_2,\cP_3 \}$ with $\cP
_1= \{\{
1\},\{2\},\{3,4\} \}$, $\cP_2= \{\{1, 2\},\{3,4\} \}$,
$\cP_3= \{\{1,2,3,4\} \}$, then $\widehat{f}_n(x_0)$ is rate
adaptive estimator of $f(x_0)$ over
\[
\bigl\{\overline{N}_{p,4}(\beta,L,\cP), (\beta,p,L,\cP)\in (0,
\beta_{\mathrm{max}} ]^4\times[1,\infty]^4\times(0,
\infty)^4\times \overline{\mP}, r (\beta,p,\overline{\varnothing} )>0
\bigr\}. %
\]
In this case, the adaptive rate of convergence satisfies
\begin{eqnarray*}
\psi_n (\beta,p,\cP )&:=& \biggl(\frac{1}{n}
\biggr)^{{r_{\mathrm{max}}}/{(2r_{\mathrm{max}}+1)}},\qquad
 (\beta,p,\cP )\in \bigl\{\beta ^{(\mathrm{max})} \bigr\}
\times \bigl\{p^{(\mathrm{max})} \bigr\}\times \{\cP _1,\cP
_2 \},\\
 r_{\mathrm{max}}&:=&\frac{\beta_{\mathrm{max}}}{2}. %
\end{eqnarray*}
Thus, in the considered example the aforementioned set contains two elements.

Finally, let us note that there is a ``$\ln$-price'' to pay for
adaptation with respect to the structure of independence even if the
smoothness parameters $\beta$, $L$ and $p$ are known. This result
follows from the bound (\ref{eq:adaptivelowereq1}) established in the
proof of Theorem~\ref{theo:adaptiveoptimality}.

\section{\texorpdfstring{Discussion: Comparison with the
global method in Lepski~\cite{lepski:supnormlossdensityestimation}}
{Discussion: Comparison with the global method in
Lepski [29]}}
\label{sec:discussion}

The latter paper deals with the rate optimal adaptive estimation of a
probability density under sup-norm loss. It is obvious that the
estimator constructed in Lepski \cite{lepski:supnormlossdensityestimation} is fully data-driven and can be
also used in pointwise estimation. However, this estimator is neither
minimax nor optimally minimax adaptive when pointwise estimation is
considered. Below, we discuss this issue in detail.

\textit{Oracle approach}. Obviously, the use of a local method allows
to control better the error of approximation since $\cB_{(h,\cP)}(x_0)$
is smaller than $\sup_{x\in\bR^d}\cB_{(h,\cP)}(x)$. Moreover, our local
method controls better the stochastic error since $\ln\delta(h,\cP)$ is
smaller than $\ln(n)$. The latter fact is explained by the use of
different constructions of the selection rule. First, it concerns the
choice of the regularization parameter $h$. Whereas Lepski \cite{lepski:supnormlossdensityestimation} uses kernel convolution, we use
the ``operation'' $\vee$ on the set of bandwidth parameters. Next, in
pointwise estimation, we select the parameter $(h,\cP)$ from very
special set whose construction is new. It is important to emphasize
that the consideration of the parameter set used in Lepski \cite{lepski:supnormlossdensityestimation} is too
``rough'' in order to bring
an optimal pointwise adaptive estimator. Both reasons required the
introduction of novel technical arguments for pointwise estimation with
respect to those in Lepski \cite{lepski:supnormlossdensityestimation}
for estimation under sup-norm loss; see the definition of our selection
rule in Section~\ref{sec:selectionrule}, and the proofs of Proposition~\ref{prop:empiricalupperbound1}, Lemma~\ref{lem:empiricalupperbound}
and Theorem~\ref{theo:oracleinequality} in the next section. Note,
however, that the adaptation to eventual independence structure in both
papers has rest upon the same methodology.

The following example illustrates clearly how the quality of estimation
provided by Lepski's estimator can be significantly improved by
application of our local method.

\begin{example} Considering the problem described in Example~\ref
{ex:densityestimation}, we compare both methods.
\begin{itemize}
\item\textit{Local method.} We obtain from our local oracle
inequality that
\[
\bigl(\bE_f^{(n)}\bigl\llvert \widehat{f}_n(x_0)-f_d(x_0)
\bigr\rrvert ^q \bigr)^{{1}/{q}}\leq\bigl(\alpha_1
\sqrt{4\ln (4 )}+\alpha _2\bigr)n^{-{1}/{2}},\qquad
\alpha_1,\alpha_2>0. %
\]

\item\textit{Global method.} The best quality of estimation provided
by Theorem~1 in Lepski \cite{lepski:supnormlossdensityestimation} is
\[
\bigl(\bE_f^{(n)}\bigl\llvert \widetilde{f}_n(x_0)-f_d(x_0)
\bigr\rrvert ^q \bigr)^{{1}/{q}}\leq(2C_1+
C_2) \bigl(n/\ln(n)\bigr)^{-{1}/{3}},\qquad C_1,C_2>0.
\]
\end{itemize}
\end{example}

It is also important to emphasize that our Theorem~\ref{theo:oracleinequality} presents other
advantages with respect to that in Lepski
\cite{lepski:supnormlossdensityestimation}.

(a) We derive our oracle-type inequality over the functional class $\bF
_d [\mathbf{f},\overline{\mP} ]$ which contains the class
$\bF
_d [\mathbf{f} ]$ used in Lepski \cite{lepski:supnormlossdensityestimation} that allows to obtain upper
bounds under more general assumptions. For instance, if $\overline{\mP
}=\{\{1,\ldots,d\}\}$, we do not need that all marginals are uniformly
bounded, that is not true when we use Theorem~1 in Lepski \cite{lepski:supnormlossdensityestimation}; see our Corollary~\ref
{cor:oracleinequality} above.

(b) The oracle-type inequality for sup-norm risk cannot be used in
general for other type of loss functions. Contrary to this, the
pointwise risk can be integrated that allows to obtain the results
under $\bL_p$-loss; see, for example, Lepski, Mammen and Spokoiny
\cite{lepskietal:Lplosswhitenoisemodel} and Goldenshluger and Lepski \cite{G-L:density-Lploss-bis}. In this context, the establishing of local
oracle inequality with the term $\ln\delta(h,\cP)$ instead of $\ln(n)$
is crucial.

\textit{Minimax adaptive estimation}. Comparing the minimax rate of
convergence defined by (\ref{eq:minimaxrate}), we find a price to pay
for adaptation in the pointwise setting. This does not exist in the
estimation under sup-norm loss. Note nevertheless that this price to
pay for adaptation is not unavoidable for all values of nuisance
parameter $ (\beta,p,L,\cP )$. This explains the necessity of the
introduction of the optimality criterion presented in Section~\ref{sec:optimality}.

Let us also compare our results with those obtained in Lepski \cite{lepski:supnormlossdensityestimation}.

\begin{example} Consider that $\overline{\mP}$ still contains the
elements $\cP_1$ and $\cP_2$ defined in Example~\ref
{ex:densityestimation} and that $d=2$. Put $\beta_{\mathrm{max}}=1$.
\begin{itemize}
\item\textit{Local method.} In view of our results, our estimator
$\widehat{f}_n$ achieves the following minimax rate of convergence:
\[
\inf_{\widehat{f}}\sup_{f\in\overline{N}_{\infty,2} (\beta
^{(\mathrm{max})},L,\cP_2 )} \bigl(
\bE_f^{(n)}\bigl\llvert \widehat {f}(x_0)-f(x_0)
\bigr\rrvert ^q \bigr)^{{1}/{q}}\asymp n^{-1/3},
\]
where infimum is taken over all possible estimators.

\item\textit{Global method.} In view of the results in Lepski \cite{lepski:supnormlossdensityestimation}, the estimator $\widetilde{f}_n$
proposed in the latter paper achieves the following minimax rate of convergence:
\[
\inf_{\widehat{f}}\sup_{f\in\overline{N}_{\infty,2} (\beta
^{(\mathrm{max})},L,\cP_2 )} \bigl(
\bE^{(n)}_f \|\widehat{f}-f \| ^q_\infty
\bigr)^{{1}/{q}}\asymp\bigl(n/\ln(n)\bigr)^{-1/3}, %
\]
where infimum is taken over all possible estimators.
\end{itemize}
\end{example}

Thus, the application of the procedure from Lepski \cite{lepski:supnormlossdensityestimation} for pointwise adaptive estimation
leads to the logarithmic loss of accuracy everywhere, while our
estimator is rate optimal for some values of nuisance parameter.

\section{Proofs of main results}
\label{proofs}
The main technical tools used in the derivation of pointwise oracle
inequality given in Theorem~\ref{theo:oracleinequality} are uniform
bounds of empirical processes. We start this section with presenting of
corresponding results those proof are postponed to the \hyperref[app]{Appendix}. In
particular, we provide with the explicit expression of the constants
$\lambda_s^{(q)}[\mathbf{K},\mz], q\geq1$, used in the
selection rule (\ref{eq:selectionrule1})--(\ref{eq:selectionrule2}).
Our considerations here are mostly based on the results recently
developed in Lepski \cite{lepski:upperfunctions}.

\subsection{Constants \texorpdfstring{$\lambda_s^{(q)}[\mathbf{K},\mz]$}{$lambda_s^{(q)}[\mathbf{K},\mathfrak{z}]$}}
\label{constants}
Set for any $s\in\bN^*, q\geq1$, $\lambda_s^{(q)}[\mathbf{K},\mz
]:= \{3q+sq[1\vee\mz](1+1/\underline{\tau}) \}^{
{1}/{2}}\lambda_s^{(q)}$, where $\underline{\tau}:=\inf_{I\in\cI
_d}\tau
(\llvert I\rrvert )>0$,
\begin{eqnarray*}
\lambda_s^{(q)}:=\lambda_s^{(q)}[
\mathbf{K}]= \biggl\{ \biggl(10se^{s}+\frac
{10seL_{\mathbf{K}}}{\llVert \mathbf{K}\rrVert _{\infty}} \biggr)\vee
(48e) \biggr\} \Bigl[\sqrt{7}+7\sqrt{(1+q)\llVert \mathbf{K}\rrVert
_{\infty}^s} \Bigr]C_{s,1}^{(q)}\llVert
\mathbf{K}\rrVert _{\infty}^s
\end{eqnarray*}
and $C_{s,1}^{(q)}:=[144s\delta_*^{-2}+5q+3+36C_s]\vee1$.

Here, $\delta_*$ is the smallest solution of the equation $8\pi
^2\delta
 (1+[\ln\delta]^2 )=1$ and
\begin{eqnarray*}
C_s&:=&s\sup_{\delta>\delta_*}\frac{1}{\delta^{2}} \biggl[1+\ln {
\biggl(\frac
{9216(s+1)\delta^{2}}{[s^*(\delta)]^{2}} \biggr)} \biggr]_+ +s\sup_{\delta
>\delta_*}
\frac{1}{\delta^{2}} \biggl[1+\ln{ \biggl(\frac
{9216(s+1)\delta
}{s^*(\delta)} \biggr)} \biggr]_+,\\
s^*(\delta)&:=&\frac{(6/\pi
^2)}{1+[\ln
\delta]^2}. %
\end{eqnarray*}

\subsection{Pointwise uniform bounds of kernel-type empirical
processes} Let $s\in\bN^*, s\leq d$, and let $Y_i=(Y_{i,1},\ldots,Y_{i,s}), i\in\bN^*$, be a sequence of $\bR^s$-valued i.i.d. random
vectors defined on a complete probability space $ (\Omega,\mA,\textsf{P} )$ and having the density $g$ with respect to the
Lebesgue measure. Later on $\bP_{g}^{(n)}$ denotes the probability law
of $Y^{(n)}:= (Y_1,\ldots,Y_n )$ and $\bE_{g}^{(n)}$ is the
mathematical expectation with respect to $\bP_{g}^{(n)}$. Assume that
$\llVert g\rrVert _{\infty}\leq\mathbf{g}$ where $\mathbf{g}>0$ is a
given number.

Set $a_s^{(q)}:= (2\sqrt{1+q} [1\vee\lambda_s^{(q)}
]
)^{-2}$ and
\begin{eqnarray*}
\cH_n^{(s)}:=\prod_{i=1}^{s}
\bigl[h_i^{(\mathrm{min})}(n),h_i^{(\mathrm{max})}(n)
\bigr]\subseteq \biggl[\frac{1}{n},1 \biggr]^s,\qquad
\mH_s^{(q)}(n):= \bigl\{h\in\cH_n^{(s)}
\dvt nV_h\geq \bigl[a_s^{(q)}
\bigr]^{-1}\ln(n) \bigr\}.
\end{eqnarray*}

For any $h\in\cH_n^{(s)}$, $y_0\in\bR^s$ and $u\geq1$ set also
\begin{eqnarray*}
K(y)&:=&\prod_{i=1}^s
\mathbf{K}(y_i),\qquad V_h:=\prod
_{i=1}^{s}{h_i}, \qquad K_{h}(y):=V_{h}^{-1}
\prod_{i=1}^s\mathbf {K}(y_i/h_i)\qquad
\forall y\in\bR^s,
\\
G_h(y_0)&:=&1\vee \biggl[\int_{\bR^s}{
\bigl\llvert K_h(y-y_0)\bigr\rrvert g(y)\,\rd y}
\biggr], \qquad\widetilde{G}_h(y_0):=1\vee
\Biggl[n^{-1}\sum_{i=1}^{n}\bigl
\llvert K_h(Y_i-y_0)\bigr\rrvert
\Biggr],
\\
\cU_h^{(u)}(y_0)&:=&\sqrt{
\frac{[G_h(y_0)]^2}{nV_h} \biggl\{1\vee\ln \biggl(\frac{V_{h^{(\mathrm{max})}}}{V_h} \biggr)+u \biggr\}}.
\end{eqnarray*}

For a given $y_0\in\bR^s$ consider the empirical processes
\begin{eqnarray*}
\xi_h^{(n)}(y_0)&:=&n^{-1}\sum
_{i=1}^{n} \bigl[K_{h}
(Y_i-y_0 )-\bE_{g}^{(n)} \bigl\{
K_h (Y_i-y_0 ) \bigr\} \bigr],\qquad h\in
\cH_n^{(s)},
\\
\overline{\xi}_h^{(n)}(y_0)&:=&n^{-1}
\sum_{i=1}^{n} \bigl[\bigl\llvert
K_{h} (Y_i-y_0 )\bigr\rrvert -
\bE_{g}^{(n)} \bigl\{ \bigl\llvert K_h
(Y_i-y_0 )\bigr\rrvert \bigr\} \bigr],\qquad h\in
\cH_n^{(s)}.
\end{eqnarray*}

\begin{proposition}
\label{prop:empiricalupperbound1}
For all $q\geq1$, all integer $n\geq3$ and all number $u$ satisfying
$1\leq u\leq q\ln(n)$
\begin{eqnarray*}
\mathrm{(i)}&& \quad \bE_g^{(n)} \Bigl\{\sup_{h\in\mH_s^{(q)}(n)}
\bigl[\bigl\llvert \xi _h^{(n)}(y_0)\bigr
\rrvert -\lambda_s^{(q)}\cU_h^{(u)}(y_0)
\bigr]_+ \Bigr\}^q \leq C_s^{(q)}(\mathbf{K},
\mathbf{g}) [nV_{h^{(\mathrm{max})}} ]^{-
{q}/{2}}\mathrm{e}^{-u};
\\
\mathrm{(ii)}&&\quad \bE_g^{(n)} \biggl\{\sup_{h\in\mH_s^{(q)}(n)}
\biggl[\bigl\llvert \overline{\xi}_h^{(n)}(y_0)
\bigr\rrvert -\frac{1}{2}G_h(y_0) \biggr]_+
\biggr\}^q \leq C_s^{(q)}(\mathbf{K},\mathbf{g})
[nV_{h^{(\mathrm{max})}} ]^{-
{q}/{2}}\mathrm{e}^{-u};
\\
\mathrm{(iii)}&&\quad \Bigl(\bE_g^{(n)} \Bigl\{\sup
_{h\in\mH_s^{(q)}(n)} \bigl[G_h(y_0)-2
\widetilde{G}_h(y_0) \bigr]_+ \Bigr\}^q
\Bigr)^{
{1}/{q}}\leq2 \bigl[C_s^{(q)}(\mathbf{K},
\mathbf{g}) \bigr]^{
{1}/{q}} [nV_{h^{(\mathrm{max})}} ]^{-{1}/{2}}\mathrm{e}^{-u/q}.
\end{eqnarray*}
\end{proposition}

The expression of the constant $C_s^{(q)}(\mathbf{K},\mathbf{g})$ is
given in the proof of the proposition.

\subsection{Oracle-type inequality}
\label{oracle-inequality-proof}
\subsubsection{Auxiliary result}
For $I\in\cI_d$ and $h\in(0,1]^d$ set
\begin{eqnarray*}
b_{h_{I}}(x_{0,I})&:=&\int_{\bR^{\llvert I\rrvert }}K_{h_{I}}^{(I)}
(x_{I}-x_{0,I} )f_{I}(x_{I})\,\rd
x_{I}, \qquad\xi_{h_I}^{(n)}(x_{0,I}):=
\widehat {f}_{h_I}^{(n)}(x_{0,I})-b_{h_{I}}(x_{0,I});
\\
G_{h_I}(x_{0,I})&:=&1\vee \biggl[\int_{\bR^{\llvert I\rrvert }}{
\bigl\llvert K_{h_I}^{(I)}(x_I-x_{0,I})
\bigr\rrvert f(x_I)\,\rd x_I} \biggr],\\
G(x_0)&:=&\sup_{(h,\cP)\in\overline{\mH}[ \overline{\mP} ]}\sup_{(\eta,\cP')\in\overline{\mH}[ \overline{\mP} ]}
\sup_{I\in\cP\circ
\cP
'}G_{h_I\vee\eta_I}(x_{0,I}).
\end{eqnarray*}
For any $(h,\cP)\in(0,1]^d\times\mP$ put
\[
\cU_{(h,\cP)}(x_0):=\sqrt{\frac{[ G(x_{0}) ]^2 \{1\vee\ln
\delta
(h,\cP) \}}{nV(h,\cP)}}. %
\]
Define also $\overline{\mathbf{f}}_n(x_0):=12\lambda d^3 (2\max
\{\overline{G}_n(x_0),1\vee\mathbf{f}\llVert \mathbf{K}\rrVert
_1^d
\} )^{d^2}$ and
\[
\xi_{n}(x_0):=\sup_{(h,\cP)\in\overline{\mH}[ \overline{\mP}
]}\sup
_{(\eta,\cP')\in\overline{\mH}[ \overline{\mP} ]}\sup_{I\in\cP
\circ\cP
'} \bigl[\bigl\llvert
\xi_{h_I\vee\eta_I}^{(n)}(x_{0,I})\bigr\rrvert -\lambda \bigl\{
\cU_{(h,\cP)}(x_0)+\cU_{(\eta,\cP')}(x_0)
\bigr\} \bigr]_+. %
\]

\begin{lemma}
\label{lem:empiricalupperbound}
Set $\mathbf{f}>0$. For any $q\geq1$ there exist constants $\mathbf
{c}_i:=\mathbf{c}_i(2q,d,\mathbf{K},\mathbf{f},\mz), i=1,2,3,4$, such
that $\forall n\geq3$, $\forall f\in\bF_d [\mathbf
{f},\overline{\mP
}  ]$, $\forall(h,\cP)\in\overline{\mH}[ \overline{\mP} ]$,
$\cP\in\mP(f)$,
\begin{eqnarray*}
\mathrm{(i)}&&\quad \bigl(\bE_f^{(n)}\bigl\llvert \xi_n(x_0)
\bigr\rrvert ^{2q} \bigr)^{1/2q}\leq \mathbf{c}_1
[nV_{\mathrm{max}} ]^{-{1}/{2}}; \\
\mathrm{(ii)}&&\quad \bigl(\bE_f^{(n)}
\bigl[G(x_0)-\overline{G}_n(x_0) \bigr]_+^{2q}
\bigr)^{1/2q}\leq\mathbf{c}_2 [nV_{\mathrm{max}}
]^{-{1}/{2}};
\\
\mathrm{(iii)}&&\quad \bigl(\bE_f^{(n)}\bigl\llvert \overline{
\mathbf{f}}_n(x_0)\bigr\rrvert ^{2q}
\bigr)^{1/2q}\leq\mathbf{c}_3;\\
\mathrm{(iv)}&&\quad \bigl(
\bE_f^{(n)}\bigl\llvert \widehat{\cU}_{(h,\cP)}(x_{0})
\bigr\rrvert ^{2q} \bigr)^{1/2q}\leq\mathbf{c}_4
\cU_{(h,\cP)}(x_0).
\end{eqnarray*}
\end{lemma}

\subsubsection{Proof of Theorem \texorpdfstring{\protect\ref{theo:oracleinequality}}{1}}
We divide the proof into several steps.

(1) Let $(h,\cP)\in\overline{\mH}[ \overline
{\mP}
]$, $\cP\in\mP(f)$, be fixed. By the triangle inequality, we have
%
\begin{eqnarray}
\label{eq:pointwiseloss1} \bigl\llvert \widehat{f}_n(x_0)-f(x_0)
\bigr\rrvert &\leq& \bigl\llvert \widehat {f}_{(\widehat{h},\widehat{\cP})}^{(n)}(x_0)-
\widehat{f}_{(h,\cP
),(\widehat{h},\widehat{\cP})}^{(n)}(x_0)\bigr\rrvert +\bigl
\llvert \widehat {f}_{(h,\cP),(\widehat{h},\widehat{\cP})}^{(n)}(x_0)-\widehat
{f}_{(h,\cP
)}^{(n)}(x_0)\bigr\rrvert\nonumber\\
&&{} +\bigl\llvert
\widehat{f}_{(h,\cP
)}^{(n)}(x_0)-f(x_0)
\bigr\rrvert
\\
&\leq& 2 \bigl[\widehat{\Delta}_{(h,\cP)}(x_0)+2\Lambda
_n(x_0)\widehat {\cU}_{(h,\cP)}(x_0)
\bigr]+\bigl\llvert \widehat{f}_{(h,\cP
)}^{(n)}(x_0)-f(x_0)
\bigr\rrvert.\nonumber
\end{eqnarray}
Here, we have used that $\widehat{f}_{(h,\cP),(\widehat{h},\widehat
{\cP
})}^{(n)}(x_0)=\widehat{f}_{(\widehat{h},\widehat{\cP}),(h,\cP
)}^{(n)}(x_0)$ and the definition of $(\widehat{h},\widehat{\cP})$.

In what follows, we will use the inequality: for $m\in\bN^*$ and
$a_i,b_i\in\bR, i=\overline{1,m}$,
%
\begin{equation}
\label{eq:productinequality} \Biggl\llvert \prod_{i=1}^{m}a_i-
\prod_{i=1}^{m}b_i\Biggr
\rrvert \leq m \Bigl(\sup_{i=\overline{1,m}}\max \bigl\{\llvert
a_i\rrvert ,\llvert b_i\rrvert \bigr\}
\Bigr)^{m-1}\sup_{i=\overline{1,m}}\llvert a_i-b_i
\rrvert .
\end{equation}
Here and later, we assume that the product and the supremum over empty
set are equal to one and zero, respectively.

(2)
Since $\cP\in\mP(f)$, using (\ref{eq:productinequality}) we have
%
\begin{eqnarray}
\label{eq:pointwiseloss2} \bigl\llvert \widehat{f}_{(h,\cP)}^{(n)}(x_0)-f(x_0)
\bigr\rrvert &\leq& d \Bigl(\sup_{I\in\cP}\max \bigl\{\widehat
{G}_{h_I}(x_{0,I}),\mathbf{f} \bigr\} \Bigr)^{d-1}
\sup_{I\in\cP
}\bigl\llvert \widehat{f}_{h_I}^{(n)}(x_{0,I})-f_I(x_{0,I})
\bigr\rrvert
\nonumber
\\[-8pt]
\\[-8pt]
\nonumber
&\leq& d \bigl(\max \bigl\{\overline{G}_{n}(x_{0}),
\mathbf{f} \bigr\} \bigr)^{d-1} \bigl[\cB_{(h,\cP)}(x_0)+
\xi_{n}(x_0)+2\lambda\cU _{(h,\cP
)}(x_0)
\bigr],
\end{eqnarray}
since $\overline{G}_{n}(x_{0})\geq\widehat{G}_{h_I}(x_{0,I})\geq1$ and
$\llvert \widehat{f}_{h_I}^{(n)}(x_{0,I})-f_I(x_{0,I})\rrvert \leq
\llvert \xi
_{h_{I}}^{(n)}(x_{0,I})\rrvert +\llvert b_{h_{I}}(x_{0,I})-f_I(x_{0,I})\rrvert $, $\forall I\in\cP$.

(3) Set $\overline{\mathbf{f}}_n^{(1)}:=d [
\overline{G}_{n}(x_{0})  ]^{d(d-1)}$. For any $(\eta,\cP')\in
\overline{\mH}[ \overline{\mP} ]$, we get from the inequality (\ref
{eq:productinequality})
\[
\bigl\llvert \widehat{f}_{(h,\cP),(\eta,\cP')}^{(n)}(x_{0})-
\widehat {f}_{(\eta,\cP')}^{(n)}(x_{0})\bigr\rrvert \leq
\overline{\mathbf{f}}_n^{(1)}\sup_{I'\in\cP'}
\biggl\llvert \prod_{I\in\cP
\dvt I\cap I'\neq\varnothing}\widehat{f}_{h_{I\cap I'},\eta_{I\cap
I'}}^{(n)}(x_{0,I\cap I'})-
\widehat{f}_{\eta
_{I'}}^{(n)}(x_{0,I'})\biggr\rrvert.
\]

Introduce, for all $I\in\cI_d$ and all $\eta\in(0,1]^d$,
$b_{h_I,\eta
_I}(x_{0,I}):=\int_{\bR^{\llvert I\rrvert }}K_{h_I\vee\eta
_I}^{(I)}
(u-x_{0,I} )f_I(u)\,\rd u$.

Put also $\overline{\mathbf{f}}_n^{(2)}:=d (\max \{
\overline
{G}_{n}(x_{0}),G(x_{0}) \} )^{d-1}$. For any $(\eta,\cP
')\in
\overline{\mH}[ \overline{\mP} ]$ and any $I'\in\cP'$, in view of
(\ref{eq:productinequality}),
\begin{eqnarray*}
&&\biggl\llvert \prod_{I\in\cP\dvt I\cap I'\neq\varnothing}\widehat
{f}_{h_{I\cap
I'},\eta_{I\cap I'}}^{(n)}(x_{0,I\cap I'})-\prod
_{I\in\cP\dvt I\cap
I'\neq
\varnothing}b_{h_{I\cap I'},\eta_{I\cap I'}}(x_{0,I\cap I'})\biggr\rrvert\\
&&\quad\leq
\overline{\mathbf{f}}_n^{(2)}\sup_{I\in\cP\dvt I\cap I'\neq
\varnothing}
\bigl\llvert \xi_{h_{I\cap I'}\vee\eta_{I\cap I'}}^{(n)}(x_{0,I\cap I'})\bigr\rrvert ,
\\
&&\biggl\llvert \prod_{I\in\cP\dvt I\cap I'\neq\varnothing}b_{h_{I\cap I'},\eta
_{I\cap I'}}(x_{0,I\cap I'})-b_{\eta_{I'}}(x_{0,I'})
\biggr\rrvert \leq\overline{\mathbf{f}}_n^{(2)}
\cB_{(h,\cP)}(x_0). %
\end{eqnarray*}
For the last inequality, we have used that $\cP\in\mP(f)$ and,
therefore, for any $\eta\in(0,1]^d$ and any $I'\in\cI_d$
\[
b_{\eta_{I'}}(x_{0,I'})= \int_{\bR^{\llvert I'\rrvert }}K_{\eta
_{I'}}^{(I')}
(x_{I'}-x_{0,I'} )\prod_{I\in\cP\dvt I\cap
I'\neq
\varnothing}f_{I\cap I'}(x_{I\cap I'})
\,\rd x_{I'}= \prod_{I\in\cP
\dvt I\cap
I'\neq\varnothing}b_{\eta_{I\cap I'}}(x_{0,I\cap I'}).
\]

(4) Applying the triangle inequality, we get since
$\overline{\mathbf{f}}_n^{(2)}\geq1$ and $\cU_{(h,\cP)}(x_{0})>0$, for
any $(\eta,\cP')\in\overline{\mH}[ \overline{\mP} ]$,
\begin{eqnarray*}
&&\bigl\llvert \widehat{f}_{(h,\cP),(\eta,\cP')}^{(n)}(x_{0})-
\widehat {f}_{(\eta,\cP')}^{(n)}(x_{0})\bigr\rrvert
\\
&&\quad\leq\overline{\mathbf{f}}_n^{(1)}\sup
_{I'\in\cP'} \Bigl\{ \overline {\mathbf{f}}_n^{(2)}
\sup_{I\in\cP\dvt I\cap I'\neq\varnothing}\bigl\llvert \xi _{h_{I\cap I'}\vee\eta_{I\cap I'}}^{(n)}(x_{0,I\cap I'})
\bigr\rrvert +\overline{\mathbf{f}}_n^{(2)}
\cB_{(h,\cP)}(x_0)+\bigl\llvert \xi_{\eta
_{I'}}^{(n)}(x_{0,I'})
\bigr\rrvert \Bigr\}
\\
&&\quad\leq\overline{\mathbf{f}}_n^{(1)}\overline{\mathbf
{f}}_n^{(2)}\cB _{(h,\cP)}(x_0)+2
\overline{\mathbf{f}}_n^{(1)}\overline{\mathbf
{f}}_n^{(2)}\xi_{n}(x_0)+3\lambda
\overline{\mathbf {f}}_n^{(1)}\overline {
\mathbf{f}}_n^{(2)} \bigl\{\cU_{(\eta,\cP')}(x_{0})+
\cU_{(h,\cP
)}(x_{0}) \bigr\}.
\end{eqnarray*}

Put $\widetilde{\mathbf{f}}_n^{(2)}:=d [ 2\overline{G}_{n}(x_{0})
 ]^{d-1}$ and $\cU(x_{0}):=\sup_{(\eta,\cP')\in\overline
{\mH}[
\overline{\mP} ]}\cU_{(\eta,\cP')}(x_{0})$. We obtain that
%
\begin{eqnarray}
\label{eq:pointwiseloss3} &&\widehat{\Delta}_{(h,\cP)}(x_0)\nonumber\\
&&\quad\leq 2
\overline{\mathbf {f}}_n^{(1)}\overline{\mathbf{f}}_n^{(2)}
\bigl\{\cB_{(h,\cP
)}(x_0)+\xi _{n}(x_0)
\bigr\} + 3\lambda\overline{\mathbf{f}}_n^{(1)} \bigl\{
\cU(x_{0})+\cU _{(h,\cP
)}(x_{0}) \bigr\} \bigl[
\overline{\mathbf{f}}_n^{(2)}-\widetilde {\mathbf
{f}}_n^{(2)} \bigr]_+
\nonumber
\\
&&\quad\quad{} + 3\lambda\overline{\mathbf{f}}_n^{(1)}\widetilde{
\mathbf {f}}_n^{(2)} \Bigl\{\sup_{(\eta,\cP')\in\overline{\mH}[ \overline
{\mP}
]}
\bigl[\cU_{(\eta,\cP')}(x_{0})-\widehat{\cU}_{(\eta,\cP
')}(x_{0})
\bigr]_+
\nonumber
\\[-8pt]
\\[-8pt]
\nonumber
&&\hspace*{73pt}{} + \bigl[\cU_{(h,\cP)}(x_{0})-\widehat{\cU
}_{(h,\cP
)}(x_{0}) \bigr]_+ \Bigr\};
\\
&&\widehat{\Delta}_{(h,\cP)}(x_0)\nonumber\\
&&\quad\leq\overline{\mathbf
{f}}_n(x_0) \bigl\{ \cB_{(h,\cP)}(x_0)+
\xi_n(x_0)+ \bigl[G(x_0)-\overline
{G}_n(x_0) \bigr]_+ \bigr\},\nonumber
\end{eqnarray}
where $\overline{\mathbf{f}}_n(x_0):=12\lambda d^3 (2\max
\{
\overline{G}_n(x_0),1\vee\mathbf{f}\llVert \mathbf{K}\rrVert
_1^d \}
 )^{d^2}$, since $\lambda\wedge\llVert \mathbf{K}\rrVert
_1\geq1$,
\[
\cU_{(\eta,\cP')}(x_0)\leq \bigl(1\vee\mathbf{f}\llVert \mathbf
{K}\rrVert _1^d \bigr)\sqrt{\frac{1\vee\ln\delta(h,\cP)}{nV(h,\cP)}}\leq 1
\vee \mathbf{f}\llVert \mathbf{K}\rrVert _1^d\qquad \forall
\bigl(\eta,\cP'\bigr)\in \overline{\mH}[ \overline{\mP} ],
\]
and $[a^m-b^m]_+\leq m(\max\{a,b\})^{m-1}[a-b]_+$, $\forall a,b>0$,
$\forall m\in\bN^*$.

(5) Finally, we deduce from (\ref{eq:pointwiseloss1}),
(\ref{eq:pointwiseloss2}) and (\ref{eq:pointwiseloss3}), using again
$\lambda\wedge\llVert \mathbf{K}\rrVert _1\geq1$, that
%
\begin{eqnarray}
\label{eq:pointwiseloss4} &&\bigl\llvert \widehat{f}_n(x_0)-f(x_0)
\bigr\rrvert
\nonumber
\\[-8pt]
\\[-8pt]
\nonumber
&&\quad\leq3\overline{\mathbf {f}}_n(x_0) \bigl\{
\cB_{(h,\cP)}(x_0)+\cU_{(h,\cP)}(x_0)+
\widehat {\cU }_{(h,\cP)}(x_{0})+\xi_n(x_0)+
\bigl[G(x_0)-\overline {G}_n(x_0)
\bigr]_+ \bigr\}.\qquad\quad
\end{eqnarray}

By the Cauchy--Schwarz inequality
\begin{eqnarray*}
&&\bigl(\bE_f^{(n)}\bigl\llvert \widehat{f}_n(x_0)-f(x_0)
\bigr\rrvert ^q \bigr)^{{1}/{q}} \\
&&\quad\leq 3 \bigl(
\bE_f^{(n)}\bigl\llvert \overline{\mathbf{f}}_n(x_0)
\bigr\rrvert ^{2q} \bigr)^{{1}/{(2q)}} \bigl[\cB_{(h,\cP)}(x_0)+
\cU_{(h,\cP
)}(x_0)+ \bigl(\bE_f^{(n)}
\bigl\llvert \widehat{\cU}_{(h,\cP
)}(x_{0})\bigr\rrvert
^{2q} \bigr)^{{1}/{(2q)}}
\\
&&\hspace*{94pt}\quad\quad{}+ \bigl(\bE_f^{(n)}\bigl\llvert \xi_n(x_0)
\bigr\rrvert ^{2q} \bigr)^{{1}/{(2q)}}+ \bigl(\bE_f^{(n)}
\bigl[G(x_0)-\overline{G}_n(x_0)
\bigr]_+^{2q} \bigr)^{{1}/{(2q)}}\bigr].
\end{eqnarray*}
Applying Lemma~\ref{lem:empiricalupperbound},
\[
\bigl(\bE_f^{(n)}\bigl\llvert \widehat{f}_n(x_0)-f(x_0)
\bigr\rrvert ^q \bigr)^{1/q}\leq3\mathbf{c}_3
\bigl[\cB_{(h,\cP)}(x_0)+(1+\mathbf {c}_4)\cU
_{(h,\cP)}(x_0)+(\mathbf{c}_1 +
\mathbf{c}_2) [nV_{\mathrm{max}} ]^{-{1}/{2}} \bigr], %
\]
and we come to the assertion of Theorem~\ref{theo:oracleinequality}
with $\alpha_1=3\mathbf{c}_3(1+\mathbf{c}_4)(1\vee\mathbf{f}\llVert
\mathbf{K}\rrVert _1^d)$ and $\alpha_2=3\mathbf{c}_3(\mathbf
{c}_1+\mathbf{c}_2)$.

\subsection{Lower bound for minimax estimation}
\subsubsection{Auxiliary result}
The result formulated in Lemma~\ref{lem:minimaxlowerbound} below is a
direct consequence of the general bound obtain in Kerkyacharian, Lepski
and Picard \cite{kerksparcecase}, Proposition~7.

Let $ (\beta,p,\cP )\in (0,\infty )^d\times
[1,\infty
]^d\times\mP$ and $L\in(0,\infty)^d$ be fixed.

\begin{lemma}
\label{lem:minimaxlowerbound}Suppose that there exists $ \{f_0,
f_1 \}\subset N_{p,d}^*(\beta,L,\cP)$ such that $\bP_{f_1}^{(n)}$
is absolutely continuous with respect to $\bP_{f_0}^{(n)}$ and
%
\begin{eqnarray}
\label{eq:lowerboundassumption1} \bigl\llvert f_1(x_0)-f_0(x_0)
\bigr\rrvert &\geq& s_n(\beta,p,\cP);
\\
\label{eq:lowerboundassumption2} \limsup_{n\rightarrow+\infty}\bE_{f_0}^{(n)}
\biggl[\frac{\rd\bP
_{f_1} ^{(n)}}{\rd\bP_{f_0}^{(n)}} \bigl(X^{(n)} \bigr)-1 \biggr]^2
&\leq& C<\infty.
\end{eqnarray}
Then, for all $q\geq1$,
\begin{eqnarray*}
&&\liminf_{n\rightarrow+\infty} \Bigl\{s_n^{-1} (\beta,p,
\cP )\inf_{\widetilde{f}_n}\cR_n^{(q)} \bigl[
\widetilde {f}_n,N_{p,d}^*(\beta,L,\cP) \bigr] \Bigr\}\\
&&\qquad\geq
\frac{1}{2} \bigl(1-\sqrt {C/(C+4)} \bigr),
\end{eqnarray*}
where infimum is taken over all possible estimators.
\end{lemma}

\subsubsection{Proof of Proposition \texorpdfstring{\protect\ref{theo:minimaxlowerbound}}{1}} 
Set $\cN(x):=\prod_{i=1}^d\sqrt{2\pi}^{
-1}\exp (-x_i^2/2 )$ and let $f_0(x):=\sigma^{-1}\cN
(x/\sigma
)$. It is easily seen that one can find $\sigma>0$ such that
%
\begin{eqnarray}
\label{eq:minimaxlowerbound0} f_0&\in& N_{p,d}^*(\beta,\underline{L}/2,
\cP)\subseteq N_{p,d}^*(\beta,L,\cP),
\nonumber
\\[-8pt]
\\[-8pt]
\nonumber
 \underline{L}_i&:=&2
\wedge L_i,\qquad  i=\overline{1,d}.
\end{eqnarray}

Let $I= \{i_1,\ldots,i_m \}\in\cP$ be such that
$r:=r
(\beta,p,\cP )=\gamma_I(\beta,p)$ and $g\dvtx\bR\rightarrow\bR
$ such
that $\operatorname{supp}(g)\subseteq(-1/2,1/2)$, $g\in\bigcap_{i\in I}N_{p_i,1}(\beta
_i,1/2)$, $\int g=0$, and $\llvert g(0)\rrvert =\llVert g\rrVert
_{\infty}$.
Define
\[
G(x_I)=A_n\prod_{l=1}^mg
\biggl(\frac{x_{i_l}-x_{0,i_l}}{\delta
_{l,n}} \biggr), %
\]
where $A_n$, $\delta_{l,n} \rightarrow0, l=\overline{1,m}$, if
$n\rightarrow\infty$, will be chosen later. Note that $G\in\bN
_{p_I,\llvert I\rrvert } (\beta_I,\underline{L}_I/2 )$ if
%
\begin{equation}
\label{eq:minimaxlowerbound1} A_n\delta_{l,n}^{-\beta_{i_l}} \Biggl(\prod
_{j=1}^m\delta _{j,n}
\Biggr)^{1/p_{i_l}}\leq\frac{\underline{L}_{i_l}}{c_l},\qquad l=\overline {1,m},
c_l=\llVert g\rrVert _{p_{i_l}}^{m-1}.
\end{equation}

Introduce
%
\begin{equation}
\label{eq:f1} f_1(x)=\prod_{i\notin I}
\bigl\{ \bigl[2\pi\sigma^2 \bigr]^{-1}\exp
\bigl(-x_i^2/2\sigma^2 \bigr) \bigr\} \biggl
\{\prod_{i\in I} \bigl[2\pi \sigma ^2
\bigr]^{-1}\exp \bigl(-x_i^2/2
\sigma^2 \bigr)+G(x_I) \biggr\}.
\end{equation}
It is obvious that there exists $A_0>0$ such that if $A_n\leq A_0$ then
$f_1(x)>0$ for any $x\in\bR^d$.
Note also that the condition $\int g=0$ implies that $\int f_1=1$. We
conclude that $f_1$ is a probability density. Furthermore, assumptions
(\ref{eq:minimaxlowerbound0})--(\ref{eq:minimaxlowerbound1}) and the
definition of $f_0$ allow us to assert that $f_1\in N_{p,d}^*(\beta,L,\cP)$.
We remark that
\[
\bigl\llvert f_1(x_0)-f_0(x_0)
\bigr\rrvert =c_1^{*} A_n,\qquad
c_1^{*}:= (\sigma \sqrt{2\pi} )^{m-d}\bigl
\llvert g(0)\bigr\rrvert ^{m}\prod_{i\notin
I}
\exp \bigl(-x_{0,i}^2/2\sigma^2 \bigr).
\]
Then Assumption (\ref{eq:lowerboundassumption1}) of Lemma~\ref
{lem:minimaxlowerbound} is fulfilled when $s_n(\beta,p,\cP)\leq
c_1^{*} A_n$.

Since $X_k, k=\overline{1,n}$, are i.i.d. random fields and $\int
g=0$ it is easily check that
\begin{eqnarray*}
\bE_{f_0}^{(n)} \biggl[\frac{\rd\bP_{f_1}^{(n)}}{\rd\bP
_{f_0}^{(n)}}
\bigl(X^{(n)} \bigr) \biggr]^2 &\leq& \Biggl[1+
\frac{2}{f_{0,I}(x_{0,I})}A_n^2 \Biggl(\prod
_{j=1}^m\delta _{j,n} \Biggr)\llVert g
\rrVert _2^{2m} \Biggr]^n \\
&\leq&\exp \Biggl[
\frac{2\llVert g\rrVert
_2^{2m}}{f_{0,I}(x_{0,I})}nA_n^2 \Biggl(\prod
_{j=1}^m\delta _{j,n} \Biggr) \Biggr],
\end{eqnarray*}
for $n$ large enough. Here, we have used that $\operatorname{supp}(G)\subseteq\Pi
_n:=\prod_{l=1}^m[x_{0,i_l}-\delta_{l,n}/2,x_{0,i_l}+\delta_{l,n}/2]$
and that $\inf_{x_I\in\Pi_n}f_{0,I}(x_I)\geq f_{0,I}(x_{0,I})/2$ for
$n$ large enough.

Since $\bE_{f_0}^{(n)} [\frac{\rd\bP_{f_1}^{(n)}}{\rd\bP
_{f_0}^{(n)}} (X^{(n)} ) ]=1$, Assumption~(\ref
{eq:lowerboundassumption2}) of Lemma~\ref{lem:minimaxlowerbound} is
fulfilled if
\[
\exp \Biggl[\frac{2\llVert g\rrVert
_2^{2m}}{f_{0,I}(x_{0,I})}nA_n^2 \Biggl(\prod
_{j=1}^m\delta_{j,n} \Biggr)
\Biggr]-1\leq C. %
\]
The latter inequality holds if
%
\begin{equation}
nA_n^2 \Biggl(\prod_{j=1}^m
\delta_{j,n} \Biggr)\leq t^2,\qquad t:=\sqrt {
\bigl[c_2^{*}\bigr]^{-1}\ln (C+1 )},\qquad
c_2^{*}:=\frac{2\llVert
g\rrVert _2^{2m}}{f_{0,I}(x_{0,I})}.
\end{equation}

To finalize our proof, we study separately two cases: $r>0$ and $r\leq
0$. Note first that $r=(1-1/s_I)/(1/\beta_I)$, where
\[
\frac{1}{s_I}:=\sum_{i\in I}\frac{1}{\beta_{i}p_{i}},\qquad
\frac
{1}{\beta
_I}:=\sum_{i\in I}\frac{1}{\beta_{i}},
\]

(1) \textit{Case} $r>0$. Solving the system
\[
A_n\delta_{l,n}^{-\beta_{i_l}} \Biggl(\prod
_{j=1}^m\delta _{j,n}
\Biggr)^{1/p_{i_l}}=\frac{\underline{L}_{i_l}}{c_l},\qquad
 l=\overline {1,m},\qquad nA_n^2
\Biggl(\prod_{j=1}^m\delta_{j,n}
\Biggr)= t^2, %
\]
we obtain
\begin{eqnarray*}
\delta_{l,n}&=& \biggl(\frac{c_l}{\underline{L}_{i_l}} \biggr)^{{1}/
{\beta_{i_l}}} \biggl(
\frac{t^2}{n} \biggr)^{{1}/{(\beta
_{i_l}p_{i_l})}}A_n^{{1}/{\beta_{i_l}}-{2}/{(\beta
_{i_l}p_{i_l})}},\qquad
A_n=R \biggl(\frac{t^2}{n} \biggr)^{{r}/{(2r+1)}},\\
 R&=& \Biggl[
\prod_{l=1}^m \biggl(\frac{\underline{L}_{i_l}}{c_l}
\biggr)^{{1}/{(2\beta_{i_l})}} \Biggr]^
{{1}/{(1-1/s_I-1/2\beta_I)}}. %
\end{eqnarray*}
It is easily seen that $A_n, \delta_{l,n} \rightarrow0, l=\overline
{1,m}$, if $n\rightarrow\infty$ and one can choose $C=1$.

We conclude that, if $r>0$, Lemma~\ref{lem:minimaxlowerbound} is
applicable with $s_n(\beta,p,\cP)= c_1^{*}R (\frac
{t^2}{n}
)^{{r}/{(2r+1)}}$.

(2) \textit{Case} $r\leq0$. We choose $A_n\equiv A$,
where the constant $A$ satisfies $0<A<A_0$. Solving the system
\begin{eqnarray*}
A\delta_{l,n}^{-\beta_{i_l}} \Biggl(\prod
_{j=1}^m\delta_{l,n}
\Biggr)^{1/p_{i_l}}&\leq&\frac{\underline{L}_{i_l}}{c_l},\qquad
 l=\overline {1,m},\qquad nA^2
\Biggl(\prod_{j=1}^m\delta_{j,n}
\Biggr)\leq t^2, %
\\
\delta_{l,n}&\geq& \biggl(\frac{Ac_l}{\underline{L}_{i_l}}
 \biggr)^{{1}/{\beta_{i_l}}}
\Biggl(\prod_{j=1}^m\delta_{j,n}
\Biggr)^{{1}/{(p_{i_l}\beta_{i_l})}}, \qquad\prod_{j=1}^m
\delta_{j,n}\leq R_2n^{-1}, \\
 R_2&=&
\frac{\ln(C+1)}{c_2^*A^2}. %
\end{eqnarray*}
Note that one can choose $A$ such that $\max_{l=\overline{1,m}}
(\frac{Ac_l}{\underline{L}_{i_l}} )^{{1}/{\beta
_{i_l}}}\leq1$
and $C=1$.
Since $s_I\leq1$, we obtain the following solution:
\[
\delta_{l,n}= \biggl(\frac{R_2}{n} \biggr)^{{s_I}/{(p_{i_l}\beta
_{i_l})}}
\rightarrow0,\qquad l=\overline{1,m}, n\rightarrow\infty. %
\]

We conclude that, if $r\leq0$, Lemma~\ref{lem:minimaxlowerbound} is
applicable with $s_n(\beta,p,\cP)= c_1^{*}A$.

This completes the proof of Proposition~\ref{theo:minimaxlowerbound}.

\subsection{Upper bounds for minimax and adaptive minimax estimation}

The proof of Theorems \ref{theo:minimaxupperbound} and \ref
{theo:adaptiveupperbound} is based on application of Theorem~\ref
{theo:oracleinequality}. Note that in view of the embedding theorem for
anisotropic Nikolskii classes (formulated in the proof of Lemma~\ref
{lem:biasupperbound}), there exists a number $\mathbf{f}:=\mathbf
{f}(\beta,p)>0$ such that $\sup_{I\in\cP}\llVert f_I\rrVert
_{\infty}\leq
\mathbf{f}$ if $r(\beta,p,\cP)>0$ or such that $\sup_{\cP\in
\overline
{\mP}^*}\sup_{I\in\cP}\llVert f_I\rrVert _{\infty}\leq\mathbf
{f}$ if
$r (\beta,p,\overline{\varnothing} )>0$. It makes possible the
application of Theorem~\ref{theo:oracleinequality}.

\subsubsection{Auxiliary result}
The result formulated in Lemma~\ref{lem:biasupperbound} below is a
consequence of Theorem~6.9 in Nikolskii \cite{nikolskii}.

Let $l\geq2$ be a fixed integer and $\overline{\mP}\subseteq\mP$
be a
fixed set of partitions of $\{1,\ldots,d\}$. Let $f\in\overline
{N}_{p,d}(\beta,L,\cP)$, where $\beta\in(0,l]^d$, $\cP\in
\overline{\mP
}$, $p\in[0,\infty]^d$ satisfy $r (\beta,p,\cP )>0$ and
$L\in
(0,\infty)^d$.

\begin{lemma}
\label{lem:biasupperbound}
There exists $\mathbf{c}:=\mathbf{c}(\mathbf{K},d,p,l,\cP)>0$ such that
\begin{eqnarray*}
\cB_{h_I,\eta_I}(x_{0,I})\leq\mathbf{c}\sum
_{i\in I}L_ih_i^{\beta
_i(I)}\qquad
\forall\cP'\in\overline{\mP}, \forall I\in\cP\circ
\cP', \forall(h,\eta)\in(0,1]^d\times[0,1]^d,
\end{eqnarray*}
where $\cB_{h_I,\eta_I}(x_{0,I})$ is defined in Section~\ref{sec:oracleinequality}, $\beta_i(I):=\kappa(I)\beta_i\kappa_i^{-1}(I)$,
$\kappa(I):=1-\sum_{k\in I} (\beta_kp_k )^{-1}$ and
$\kappa
_i(I):=1-\sum_{k\in I} (p_k^{-1}-p_i^{-1} )\beta_k^{-1}$.
\end{lemma}

The proof of this lemma is given in the \hyperref[app]{Appendix}.

\subsubsection{Proof of Theorem \texorpdfstring{\protect\ref{theo:minimaxupperbound}}{3}}
For all $I\in\cP$, consider the following system of equations:
\[
h_j^{\beta_j(I)}=h_i^{\beta_i(I)}=\sqrt{
\frac{1}{nV_{h_I}}}, \qquad i,j\in I, %
\]
and let $\mathbf{h}_I$ denotes its solution. One can easily check that
%
\begin{equation}
\label{eq:minimaxbandwidth1} \mathbf{h}_i=
n^{-({\gamma_I(\beta,p)}/{(2\gamma_I(\beta,p)+1)})
({1}/{\beta_i(I)})}, \qquad i\in I, I\in\cP.
\end{equation}
Here, we have used that $1/\gamma_I(\beta,p)=\sum_{i\in I}1/\beta_i(I)$.

We note that $2^{-1}nV(\mathbf{h},\cP)\geq a^{-1}\ln(n)$ for all $n$ large enough. To get the statement of the
theorem, we will apply Theorem~\ref{theo:oracleinequality}
with $\mz=1$, $\tau(s)=1$, $s=1,\ldots,d$, $\mh_I^{(I)}=\mathbf{h}_I$ if $I\in\cP$ and $\mh_i^{(I)}=1$
if $i\in I$, $I\notin\cP$, $\overline{\mH}=\{\mathbf{h}\}$, $\overline{\mP}=\{\cP\}$.
Thus, $\overline{\mH}[\overline{\mP}]$ is non-empty for $n$ large enough and we get
%
\begin{equation}
\label{eq:minimaxupperbound1} \cR_n^{(q)} \bigl[\widehat{f}_{(\mathbf{h},\cP)}^{(n)},f
\bigr]\leq \alpha_1 (\mathbf{c}\overline{L}\vee1 ) \biggl[ \sup
_{I\in\cP
}\sum_{i\in I}
\mathbf{h}_i^{\beta_i(I)}+\sup_{I\in\cP}\sqrt {
\frac
{1}{nV_{\mathbf{h}_I}}} \biggr]+\alpha_2 \sup_{I\in\cP}
\sqrt {\frac
{1}{nV_{\mathbf{h}_I}}},
\end{equation}
where $\overline{L}:=\sup_{i=\overline{1,d}} L_i$. Here, we have used
Lemma~\ref{lem:biasupperbound} and the definition of $\cB_{(\mathbf
{h},\cP)}(x_0)$.

We deduce from (\ref{eq:minimaxbandwidth1}) and (\ref{eq:minimaxupperbound1})
\[
\cR_n^{(q)} \bigl[\widehat{f}_{(\mathbf{h},\cP)}^{(n)},f
\bigr]\leq \bigl[2\alpha_1 (\mathbf{c}\overline{L}\vee1 )+\alpha
_2 \bigr]\sup_{I\in\cP}
n^{-{\gamma_I(\beta,p)}/{(2\gamma_I(\beta,p)+1)}}= \bigl[2
\alpha_1 (\mathbf{c}\overline{L}\vee1 )+\alpha _2
\bigr]n^{-{r}/{(2r+1)}}
\]
and the assertion of Theorem~\ref{theo:minimaxupperbound} follows.

\subsubsection{Proof of Theorem \texorpdfstring{\protect\ref{theo:adaptiveupperbound}}{4}}
Set $ (\beta,p )\in (0,\beta_{\mathrm{max}} ]^d\times
[1,\infty
]^d$ such that $r (\beta,p,\overline{\varnothing} )>0$, $\cP
\in
\overline{\mP}$, $L\in(0,\infty)^d$, and $f\in\overline
{N}_{p,d}(\beta,L,\cP)$.

Let us first note the following simple fact. If $\cP'\in\overline
{\mP}$
and $J=I\cap I', I\in\cP, I'\in\cP'$, we easily prove that $\beta
_i(J)\geq\beta_i(I)\ \forall i\in J$; see, for example, Lepski \cite{lepski:supnormlossdensityestimation}, proof of Theorem~3, for more
details. Thus, in view of Lemma~\ref{lem:biasupperbound},
%
\begin{equation}
\label{eq:adaptiveupperbound1} \cB_{(h,\cP)}(x_0)\leq\mathbf{c}\sup
_{I\in\cP}\sum_{i\in
I}L_ih_i^{\beta
_i(I)}\qquad
\forall h\in(0,1]^d.
\end{equation}

Recall that $\mh_I^{(I)}$, $I\in\cI_d$, is the projection on the dyadic grid in $(0,1]^{|I|}$ of $\rh_I^{(I)}$ given in
(\ref{eq:maximalebandwidth}) and note that $2^{-1}nV_{\mh_I^{(I)}}\geq a^{-1}\ln(n)$ for $n$ large enough.
Thus, $\overline{\mH}[\overline{\mP}]$ is non-empty and one can apply Theorem \ref{theo:oracleinequality}.


If $r (\beta,p,\cP )=r_{\mathrm{max}}$, then it is obvious that
$
(\beta,p )= (\beta^{(\mathrm{max})},p^{(\mathrm{max})} )$ and that
$d(\cP
)=\overline{d}$. Thus, in view of the definition of the multibandwidths
$\mathfrak{h}_I^{(I)}$, $I\in\cP$, $\inf_{I\in\cP}V_{\mathfrak
{h}_I^{(I)}}=V_{\mathrm{max}}$. It
follows from Theorem~\ref{theo:oracleinequality} and (\ref
{eq:adaptiveupperbound1})
\[
\cR_n^{(q)} [\widehat{f}_n,f ]\leq
\alpha_1 (\mathbf {c}\overline{L}\vee1 ) \biggl[ \sup
_{I\in\cP}\sum_{i\in
I} \bigl(
\mathfrak{h} _i^{(I)} \bigr)^{\beta_{\mathrm{max}}}+\sup
_{I\in\cP}\sqrt{\frac
{1}{nV_{\mathfrak{h}
_I^{(I)}}}} \biggr]+\alpha_2
[nV_{\mathrm{max}} ]^{-{1}/{2}}, %
\]
where $\overline{L}:=\sup_{i=\overline{1,d}} L_i$. Since
$r_{\mathrm{max}}=\beta
_{\mathrm{max}}/\overline{d}$, we conclude that there exists a constant $C>0$
such that
%
\begin{equation}
\label{eq:adaptiveupperbound2} \cR_n^{(q)} [\widehat{f}_n,f ]
\leq C \bigl[\alpha_1 (\mathbf {c}\overline{L}\vee1 ) (d+1)+
\alpha_2 \bigr]n^{-{r_{\mathrm{max}}}/{(2r_{\mathrm{max}}+1)}}.
\end{equation}

If $r (\beta,p,\cP )<r_{\mathrm{max}}$ we solve, for all $I\in\cP
$, the system
\[
L_jh_j^{\beta_j(I)}=L_ih_i^{\beta_i(I)}=
\sqrt{\frac{\ln
(n)}{nV_{h_I}}},\qquad i,j\in I. %
\]
The solution is
%
\begin{eqnarray}
\label{eq:adaptiveminimaxbandwidth1} h_i&=&
L_i^{-{1}/{\beta_i(I)}} \biggl(
\frac{L(I)\ln(n)}{n} \biggr)^{
{\gamma_I(\beta,p)}/{(2\gamma_I(\beta,p)+1)}
{1}/{\beta_i(I)}},
\nonumber
\\[-8pt]
\\[-8pt]
\nonumber
 L(I)&=&\prod
_{i\in I}L_i^{{1}/{\beta_i(I)}},\qquad i\in I, I\in \cP.
\end{eqnarray}
It is easily seen that $(h,\cP)\in\mH[ \mP]$ for $n$ large enough.
Replacing $h$ by its projection $\bar{h}$ on the dyadic grid
$\overline{\mH}$, one has $(\bar{h},\cP)\in\overline{\mH}[
\overline{\mP} ]$ for $n$ large enough. We deduce from Theorem~\ref
{theo:oracleinequality} and (\ref{eq:adaptiveupperbound1})
%
\begin{equation}
\label{eq:adaptiveupperbound3} \cR_n^{(q)} [\widehat{f}_n,f ]
\leq\alpha_1 \biggl[ \mathbf {c}\sup_{I\in\cP}\sum
_{i\in I}L_i\bar{h}_i^{\beta_i(I)}+
\sup_{I\in
\cP}\sqrt{\frac{\ln(n)}{nV_{\bar{h}_I}}} \biggr]+
\alpha_2 [nV_{\mathrm{max}} ]^{-{1}/{2}}.
\end{equation}

The assertion of Theorem~\ref{theo:adaptiveupperbound} follows from
(\ref{eq:adaptiveupperbound2}), (\ref{eq:adaptiveminimaxbandwidth1})
and (\ref{eq:adaptiveupperbound3}).

\subsection{Lower bound for adaptive minimax estimation and optimal rate}
\subsubsection{Auxiliary result}
To get the assertion of Theorem~\ref{theo:adaptiveoptimality}, we use
the following lemma which is due to an oral communication with O.
Lepski. This result can be viewed as a generalization of Lemma~\ref
{lem:minimaxlowerbound}.

Let $ (\beta,p )\in (0,\beta_{\mathrm{max}} ]^d\times
[1,\infty
]^d$ such that $r (\beta,p,\overline{\varnothing} )>0$, $\cP
\in
\overline{\mP}$, $L\in(0,\infty)^d$ and $ (\beta',p'
)\in
(0,\beta_{\mathrm{max}} ]^d\times[1,\infty]^d$ such that $r (\beta
',p',\overline{\varnothing} )>0$, $\cP'\in\overline{\mP}$,
$L'\in
(0,\infty)^d$ be fixed.

\begin{lemma}
\label{lem:adaptivelowerbound}
Set $(a_n)$ and $(b_n)$ two sequences such that $a_n, b_n,
b_n/a_n\rightarrow\infty, n\rightarrow\infty$.
Suppose that exist $f_0\in N_2:=\overline{N}_{p',d}(\beta',L',\cP')$
and $f_1\in N_1:=\overline{N}_{p,d}(\beta,L,\cP)$ such that $\bP
_{f_1}^{(n)}$ is absolutely continuous with respect to $\bP
_{f_0}^{(n)}$ and
%
\begin{equation}
\label{eq:adaptivelowerbound0} \bigl\llvert f_1(x_0)-f_0(x_0)
\bigr\rrvert = a_n^{-1}; \qquad \bE_{f_0}^{(n)}
\biggl[\frac{\rd\bP_{f_1}^{(n)}}{\rd\bP
_{f_0}^{(n)}} \bigl(X^{(n)} \bigr) \biggr]^2\leq
\frac{b_n}{a_n}.
\end{equation}
Then, for any $q\geq1$,
\begin{eqnarray*}
\liminf_{n\rightarrow+\infty}\inf_{\widetilde{f}_n} \Bigl[\sup
_{f\in
N_1}\bE_f^{(n)} \bigl\{a_n
\bigl|\widetilde{f}_n(x_0)-f(x_0) \bigr| \bigr
\}^q +\sup_{f\in N_2} \bE_f^{(n)}
\bigl\{b_n \bigl|\widetilde {f}_n(x_0)-f(x_0)
\bigr| \bigr\}^q \Bigr]\geq\frac{1}{2},
\end{eqnarray*}
where infimum is taken over all possible estimators.
\end{lemma}

The proof of this lemma is given in the \hyperref[app]{Appendix}.

\subsubsection{Proof of Theorem \texorpdfstring{\protect\ref{theo:adaptiveoptimality}}{5}}

(1) Set $N_1:=\overline{N}_{p,d}(\beta,L,\cP)$,
$N_2:=\overline{N}_{p',d}(\beta',L',\cP')$, $r_1:=r(\beta,p,\cP)$ and
$r_2:=r(\beta',p',\cP')$ such that $0<r_1<r_2$. For any $\tau$ such
that $\frac{r_1}{2r_1+1}<\tau\leq\frac{r_2}{2r_2+1}$, there exists
$C(\tau)>0$ satisfying: $\forall q\geq1$,
%
\begin{eqnarray}
\label{eq:adaptivelowereq1} &&\liminf_{n\rightarrow+\infty}\inf_{\widetilde{f}_n}
\biggl[\sup_{f\in
N_1}\bE_f^{(n)} \biggl\{
\biggl(\frac{n}{\ln(n)} \biggr)^{{r_1}/{(2r_1
+1)}}\bigl |\widetilde{f}_n(x_0)-f(x_0)
\bigr| \biggr\}^q
\nonumber
\\[-8pt]
\\[-8pt]
\nonumber
&&\hspace*{22pt}\qquad{}+\sup_{f\in N_2} \bE_f^{(n)}
\bigl\{n^{\tau} \bigl|\widetilde {f}_n(x_0)-f(x_0)
\bigr| \bigr\}^q \biggr]\geq C(\tau).
\end{eqnarray}

Let us prove (\ref{eq:adaptivelowereq1}). The proof is based on Lemma~\ref{lem:adaptivelowerbound} where we put
\[
a_n:= \bigl[2C(\tau) \bigr]^{-1/q} \biggl(
\frac{n}{\ln(n)} \biggr)^{
{r_1}/{(2r_1 +1)}},\qquad b_n:= \bigl[2C(\tau)
\bigr]^{-1/q}n^{\tau}, %
\]
and the constant $C(\tau)>0$ will be specified later.

Similarly to the proof of Proposition~\ref{theo:minimaxlowerbound}, set
$\cN(x):=\prod_{i=1}^d\sqrt{2\pi}^{ -1}\exp (-x_i^2/2
)$ and
define $f_0(x):=\sigma^{-1}\cN(x/\sigma)$, where $\sigma$ is chosen in
such way that
\[
f_0\in\overline{N}_{p',d}\bigl(\beta',L',
\cP'\bigr)\cap\overline {N}_{p,d}(\beta,\underline{L}/2,
\cP). %
\]
Let also $f_1$ be given in (\ref{eq:f1}).
It is obvious that there exists a constant $A_0$ such that $f_1\in N_1$
if $A_n\leq A_0$ and
%
\begin{equation}
\label{eq:adaptiveminimaxlowerbound1} A_n\delta_{l,n}^{-\beta_{i_l}} \Biggl(\prod
_{j=1}^m\delta _{j,n}
\Biggr)^{1/p_{i_l}}\leq\frac{\underline{L}_{i_l}}{c_l},\qquad l=\overline {1,m},
c_l=\llVert g\rrVert _{p_{i_l}}^{m-1}.
\end{equation}

Assumptions of Lemma~\ref{lem:adaptivelowerbound} are, respectively,
fulfilled if
%
\begin{eqnarray}
c_1^{*} A_n&\geq& \bigl[2C(\tau)
\bigr]^{1/q} \biggl(\frac{\ln
(n)}{n} \biggr)^{{r_1}/{(2r_1+1)}},\nonumber\\
c_1^{*}&:=& (\sigma\sqrt {2\pi } )^{m-d}\bigl
\llvert g(0)\bigr\rrvert ^{m}\prod_{i\notin I}
\exp \bigl(-x_{0,i}^2/2\sigma^2 \bigr);
\\
 \exp \Biggl[\frac{2\llVert g\rrVert
_2^{2m}}{f_{0,I}(x_{0,I})}nA_n^2 \Biggl(\prod
_{l=1}^m\delta _{l,n} \Biggr)
\Biggr]&\leq&n^{\tau} \biggl(\frac{n}{\ln(n)} \biggr)
^{-{r_1}/{(2r_1+1)}}.
\nonumber
\end{eqnarray}
The latter inequality, in its turn, holds if
%
\begin{equation}
nA_n^2 \Biggl(\prod_{l=1}^m
\delta_{l,n} \Biggr)=t^2\ln(n),\qquad t:=\sqrt {
\bigl[c_2^{*}\bigr]^{-1} \biggl(\tau-
\frac{r_1}{2r_1+1} \biggr)},\qquad c_2^{*}:=\frac{2\llVert g\rrVert _2^{2m}}{f_{0,I}(x_{0,I})}.
\end{equation}

Solving the system
\[
A_n\delta_{l,n}^{-\beta_{i_l}} \Biggl(\prod
_{j=1}^m\delta _{j,n}
\Biggr)^{1/p_{i_l}} = \frac{\underline{L}_{i_l}}{c_l},\qquad
 l=\overline {1,m},\qquad
nA_n^2 \Biggl(\prod_{l=1}^m
\delta_{l,n} \Biggr)=t^2\ln(n), %
\]
we obtain
\begin{eqnarray*}
\delta_{l,n}&=& \biggl(\frac{c_l}{L_{i_l}} \biggr)^{{1}/{\beta
_{i_l}}} \biggl(
\frac{t^2\ln(n)}{n} \biggr)^{{1}/{(\beta
_{i_l}p_{i_l})}}A_n^{{1}/{\beta_{i_l}}-{2}/{(\beta
_{i_l}p_{i_l})}},\qquad
A_n=R \biggl(\frac{t^2\ln(n)}{n} \biggr)^{
{r_1}/{(2r_1+1)}},\\
  R&=& \Biggl[
\prod_{l=1}^m \biggl(\frac
{L_{i_l}}{c_l}
\biggr)^{{1}/{(2\beta_{i_l})}} \Biggr]^{{1}/{(1-1/s_I-1/2\beta_I)}}. %
\end{eqnarray*}
It is easily seen that $A_n, \delta_{l,n} \rightarrow0, l=\overline
{1,m}$, if $n\rightarrow\infty$. The choice
$
C(\tau)=\frac{1}{2} [c_1^*R (t^{{2r_1}/{(2r_1+1)}}
) ]^{q}$,
completes the proof of the inequality (\ref{eq:adaptivelowereq1}). It
follows the assertion (i) of Theorem~\ref{theo:adaptiveoptimality}.

(2) Let us recall the definition of the set $\cA
\times
\mB$, which is the set of ``nuisance'' parameters for the considered problem.
\[
\cA:= \bigl\{ (\beta,p )\in (0,\beta_{\mathrm{max}} ]^d\times [1,
\infty]^d\dvt r (\beta,p,\overline{\varnothing} )>0 \bigr\}, \qquad\mB:=
\overline{\mP}. %
\]

Let $\widetilde{\psi}_n$ be an admissible family of normalizations and
let $\widetilde{f}_n(x_0)$ be $\widetilde{\psi}_n$-adaptive
estimator. Define
\begin{eqnarray*}
\cA^{(0)} [\widetilde{\psi}/\psi ]&:=& \Bigl\{ (\beta,p )\in\cA\dvt
\lim_{n\to\infty}\Upsilon_n (\beta,p )=0 \Bigr\},
\\
\Upsilon_n (\beta,p )&:=&\inf_{\cP\in\overline{\mP
}}\Upsilon
_n (\beta,p,\cP ), \qquad\Upsilon_n (\beta,p,\cP ):=
\frac{\widetilde{\psi}_n (\beta,p,\cP )}{\psi
_n (\beta,p,\cP )},
\end{eqnarray*}
where $\psi_n$ is given in (\ref{eq:adaptiverate}).
For any $\cP\in\overline{\mP}$ put also
\[
\cA_{\cP}^{(\infty)} [\widetilde{\psi}/\psi ]:= \Bigl\{ (\beta,p )
\in\cA\dvt \lim_{n\to\infty}\Upsilon_n (\beta
_0,p_0 )\Upsilon_n (\beta,p,\cP )=\infty,
\forall (\beta _0,p_0 )\in\cA^{(0)} [
\widetilde{\varPsi}/\varPsi ] \Bigr\}. %
\]
In the slight abuse of the notation, we will use later $\psi_n(r)$
instead of $\psi_n(\beta,p,\cP)$, $r=r(\beta,p,\cP)$.

For any $(\beta_0,p_0)\in\cA^{(0)} [\widetilde{\psi}_n/\psi
_n ]$ introduce
%
\begin{equation}
\label{eq:outperformrate} \cP_0:=\arg\inf_{\cP\in\overline{\mP}}
\Upsilon_n (\beta _0,p_0,\cP ),\qquad
r_0:=r(\beta_0,p_0,
\cP_0).
\end{equation}

Let us first note that $0<r_0<r_{\mathrm{max}}$ for any $(\beta_0,p_0)\in\cA
^{(0)} [\widetilde{\psi}_n/\psi_n ]$. Indeed, if $r_0=r_{\mathrm{max}}$
then $(\beta_0,p_0)\in\cA^{(0)} [\widetilde{\psi}_n/\psi
_n ]$
contradicts to $\psi_n (r_{\mathrm{max}} )$ is a minimax rate of
convergence. Moreover, for any $r\in(r_0,r_{\mathrm{max}})$, there exists $
(\beta,p )\in\cA$ and $\cP\in\overline{\mP}$ such that
$r (\beta,p,\cP )=r$. It suffices to choose $\cP$ such that $r (\beta
^{(\mathrm{max})},p^{(\mathrm{max})},\cP )=r_{\mathrm{max}}=\beta_{\mathrm{max}}/\llvert I\rrvert $,
$I\in
\cP$, and $\beta_i=r\llvert I\rrvert $, $p_i=\infty$, $i=1,\ldots,d$.

(3) Our goal now is to prove that for any $(\beta
_0,p_0)\in\cA^{(0)} [\widetilde{\psi}_n/\psi_n ]$ we have
%
\begin{equation}
\label{eq:resultstep2} \lim_{n\to\infty}\Upsilon_n (
\beta_0,p_0 )\Upsilon _n (\beta,p,\cP )=
\infty \qquad\forall(\beta,p,\cP)\dvt r_0<r (\beta,p,\cP
)<r_{\mathrm{max}}.
\end{equation}

Set $N_0:=\overline{N}_{p_0,d}(\beta_0,L_0,\cP_0)$ and $N:=\overline
{N}_{p,d}(\beta,L,\cP)$ such that $r_0<r (\beta,p,\cP
)<r_{\mathrm{max}}$. Applying the inequality (\ref{eq:adaptivelowereq1}) with
$r_1=r_0, N_1=N_0$, $r_2=r$ and $ N_2=N$, we get for any $\tau$
satisfying $\frac{r_0}{2r_0+1}<\tau<\frac{r}{2r+1}$
%
\begin{eqnarray}
\label{eq:optimal1}&& \liminf_{n\rightarrow+\infty} \Bigl[\sup_{f\in N_0}
\bE_f^{(n)} \bigl\{ \psi_n^{-1}(r_0)
\bigl|\widetilde{f}_n(x_0)-f(x_0) \bigr| \bigr
\}^q +\sup_{f\in N} \bE_f^{(n)}
\bigl\{n^{\tau} \bigl|\widetilde {f}_n(x_0)-f(x_0)
\bigr| \bigr\}^q \Bigr]
\nonumber
\\[-8pt]
\\[-8pt]
\nonumber
&&\quad\geq C(\tau).
\end{eqnarray}
Furthermore, by definition of $\widetilde{f}_n(x_0)$ and $\widetilde
{\psi}_n$, there exist constants $M_0,M>0$ such that for all $n$ large enough
%
\begin{eqnarray}
\label{eq:optimal2} \sup_{f\in N_0}\bE_f^{(n)}
\bigl\{\widetilde{\psi}_n^{-1}(\beta _0,p_0,
\cP_0) \bigl|\widetilde{f}_n(x_0)-f(x_0)
\bigr| \bigr\}^q&\leq& M_0;
\\
\label{eq:optimal3} \sup_{f\in N} \bE_f^{(n)}
\bigl\{\widetilde{\psi}_n^{-1}(\beta,p,\cP )\bigl |
\widetilde{f}_n(x_0)-f(x_0) \bigr| \bigr
\}^q&\leq& M.
\end{eqnarray}

Note that $\lim_{n\to\infty}\frac{\widetilde{\psi}_n(\beta
_0,p_0,\cP
_0)}{\psi_n(\beta_0,p_0,\cP_0)}=0$ that follows from $(\beta
_0,p_0)\in
\cA^{(0)} [\widetilde{\psi}_n/\psi_n ]$ as well as the definition
of $\cP_0$. Thus, we obtain in view of (\ref{eq:optimal2}) that
\[
\lim_{n\to\infty}\sup_{f\in N_0}\bE_f^{(n)}
\bigl\{\psi _n^{-1}(r_0)\bigl |
\widetilde{f}_n(x_0)-f(x_0) \bigr| \bigr
\}^q=0. %
\]
It yields together with (\ref{eq:optimal1}) and (\ref{eq:optimal3}) that
%
\begin{equation}
\label{eq:optimal4} \liminf_{n\rightarrow+\infty} Mn^{\tau}\widetilde{
\psi}_n(\beta,p,\cP )\geq C(\tau).
\end{equation}

Recall that $\psi_n(r)=(\ln(n)/n)^{{r}/{(2r+1)}}$. Since $\tau
<\frac
{r}{2r+1}$ we get for some $a>0$ satisfying $\tau+a<\frac{r}{2r+1}$ that
$n^{\tau}\psi_n(r)\leq n^{-a}$ for $n$ large enough. Hence, we obtain
in view of (\ref{eq:optimal4})
%
\begin{equation}
\label{eq:optimal5} \liminf_{n\rightarrow+\infty} n^{-a}
\Upsilon_n (\beta,p,\cP ):=\liminf_{n\rightarrow+\infty}
n^{-a}\frac{\widetilde{\psi
}_n(\beta,p,\cP)}{\psi_n(\beta,p,\cP)}\geq\frac{C(\tau)}{M}.
\end{equation}
Furthermore, since $\varphi_n(\beta_0,p_0,\cP_0)$ is a minimax rate of
convergence, there exists a constant \mbox{$M_1>0$} such that
%
\begin{equation}
\label{eq:optimal6} \Upsilon_n (\beta_0,p_0
):=\frac{\widetilde{\psi
}_n(\beta
_0,p_0,\cP_0)}{\psi_n(\beta_0,p_0,\cP_0)}\geq M_1\frac{\varphi
_n(\beta
_0,p_0,\cP_0)}{\psi_n(\beta_0,p_0,\cP_0)}=M_1
\bigl[\ln(n)\bigr]^{-{r_0}/{(2r_0+1)}}
\end{equation}
for all $n$ large enough. We deduce from (\ref{eq:optimal5}) and (\ref
{eq:optimal6}) that $\lim_{n\to\infty}\Upsilon_n (\beta
_0,p_0
)\Upsilon_n (\beta,p,\cP )=\infty$.

(4) Let $ (\beta_1,p_1 )\in\cA^{(0)}
[\widetilde{\psi}_n/\psi_n ]$ and $ (\beta_2,p_2 )\in
\cA
^{(0)} [\widetilde{\psi}_n/\psi_n ]$ be arbitrary pairs of
parameters. Let also $\cP_1$ and $\cP_2$ be defined in (\ref
{eq:outperformrate}) where $ (\beta_0,p_0 )$ is replaced by
$
(\beta_1,p_1 )$ and $ (\beta_2,p_2 )$, respectively. Then
necessarily
%
\begin{equation}
\label{eq:outperformrate1} r (\beta_1,p_1,\cP_1
)=r (\beta_2,p_2,\cP_2 ).
\end{equation}
Indeed, assume that $r (\beta_1,p_1,\cP_1 )<r (\beta
_2,p_2,\cP
_2 )$. Noting that $\Upsilon_n (\beta_2,p_2 )=\Upsilon
_n (\beta_2,p_2,\cP_2 )$, in view of the definition of
$\cP_2$
we deduce from (\ref{eq:resultstep2}) with $ (\beta_1,p_1
)=
(\beta_0,p_0 )$ and $ (\beta,p,\cP )= (\beta
_2,p_2,\cP_2
)$ that
%
\begin{equation}
\label{eq:outperformrate2} \Upsilon_n (\beta_2,p_2 )
\rightarrow\infty, \qquad n\rightarrow \infty.
\end{equation}
This contradicts to $ (\beta_2,p_2 )\in\cA^{(0)}
[\widetilde
{\psi}_n/\psi_n ]$. The case $r (\beta_1,p_1,\cP_1
)>r
(\beta_2,p_2,\cP_2 )$ is traited similarly.

(5) We are now in position to prove Theorem~\ref
{theo:adaptiveoptimality}.

First, if $\cA^{(0)} [\widetilde{\psi}_n/\psi_n ]\neq
\varnothing
$, we deduce from (\ref{eq:outperformrate1}) that there exists $r_0\in
(0,r_{\mathrm{max}})$ such that
%
\begin{equation}
\label{eq:optimal7} r (\beta,p,\cP_{(\beta,p)} )=r_0\qquad \forall (
\beta,p )\in \cA^{(0)} [\widetilde{\psi}_n/\psi_n
].
\end{equation}
Here, as previously, $\cP_{(\beta,p)}:=\arg\inf_{\cP\in\overline
{\mP
}}\Upsilon_n (\beta,p,\cP )$.

Recall that, for $ (\beta,p,\cP )\in (0,+\infty
)^d\times[1,\infty]^d\times\mP$,
\[
r (\beta,p,\cP )=\inf_{I\in\cP}\gamma_I (\beta,p ),\qquad
\gamma_I (\beta,p )=\frac{1-\sum_{i\in I}
{1}/{(\beta_i p_i)}}{\sum_{i\in I}{1}/{\beta_i}},\qquad I\in\cP. %
\]
Thus, obviously
%
\begin{equation}
\label{eq:optimal8} \operatorname{dim} \bigl(\cA^{(0)} [\widetilde{
\psi}_n/\psi _n ] \bigr)\leq2d-1.
\end{equation}

Next, let $\cP^*\in\overline{\mP}$ be a partition satisfying $r
(\beta^{(\mathrm{max})},p^{(\mathrm{max})},\cP^* )=r_{\mathrm{max}}$. We deduce from (\ref
{eq:resultstep2}) that
%
\begin{equation}
\label{eq:optimal9} \cA_{\cP^*}^{(\infty)} [\widetilde{\psi}/\psi ]
\supseteq \bigl\{ (\beta,p )\in\cA\dvt r_0<r \bigl(\beta,p,\cP^*
\bigr)<r_{\mathrm{max}} \bigr\},
\end{equation}
where $r_0$ is defined in (\ref{eq:optimal7}). Thus, $\cA_{\cP
^*}^{(\infty)} [\widetilde{\psi}/\psi ]$ contains an open
set of
$\cA$ since $(\beta,p)\mapsto r (\beta,p,\cP^* )$ is
continuous. This together with (\ref{eq:optimal8}) completes the proof
of the theorem.

\begin{appendix}\label{app}

\section*{Appendix}
\setcounter{subsection}{0}
\subsection{Proof of Proposition \texorpdfstring{\protect\ref{prop:empiricalupperbound1}}{2}}
\label{sec:proofproposition1}
Our goal is to establish a uniform bound for the empirical process
$
\{\xi_h^{(n)}(y_0) \}_h$. Note that the considered family of random
fields is a particular case of the generalized empirical processes
studied in Lepski \cite{lepski:upperfunctions}. We get the assertions
of Proposition~\ref{prop:empiricalupperbound1} from the Theorem $1$ in
the latter paper since it allows us to assert that, for any $u\geq1$,
$q\geq1$ and any integer $n\geq3$
%
\renewcommand{\theequation}{\arabic{equation}}
\setcounter{equation}{57}
\begin{eqnarray}
\label{eq:empiricalupperbound1} &&\bE_g^{(n)} \Bigl\{\sup
_{h\in\cH_n^{(s)}} \bigl[\bigl\llvert \xi _h^{(n)}(y_0)
\bigr\rrvert -\cU^{(u,q)}(n,h,y_0) \bigr]_+ \Bigr
\}^q \leq C_s^{(q)}(\mathbf{K},\mathbf{g})
[nV_{h^{(\mathrm{max})}} ]^{-{q}/{2}}\mathrm{e}^{-u},
\nonumber\\
&&\cU^{(u,q)}(n,h,y_0)
\nonumber
\\[-8pt]
\\[-8pt]
\nonumber
&&\quad:=c(\mathbf{K},s,q)\sqrt{
\frac
{G_h(y_0)}{nV_h} \biggl\{1\vee\ln \biggl(\frac
{V_{h^{(\mathrm{max})}}}{V_h} \biggr)+2\ln
\bigl(2+\ln G_h(y_0) \bigr)+u \biggr\}}
\\
&&\quad\quad{} + \frac{c(\mathbf
{K},s,q)}{nV_h} \biggl\{1\vee\ln \biggl(\frac
{V_{h^{(\mathrm{max})}}}{V_h} \biggr)+2\ln
\bigl(2+\ln G_h(y_0) \bigr)+u \biggr\}.
\nonumber
\end{eqnarray}
The constants $C_s^{(q)}(\mathbf{K},\mathbf{g})$ and $c(\mathbf
{K},s,q)$ are given later.

Thus, we only have to check the Assumptions of Theorem $1$ in Lepski
\cite{lepski:upperfunctions} and to match the notation used in the
present paper and in the latter one. We divide this proof into several steps.

(1) For our case, we first consider that $p=1$,
$m=s+1$, $k=s$, $\mH_1^{k}(n)=\cH_n^{(s)}$, $\mH_{k+1}^m(n)= \{
y_0 \}$, $\mathfrak{h}^{(k)}=h$ and
\begin{eqnarray*}
G_{\infty}(h)&=&V_h^{-1}\llVert \mathbf{K}\rrVert
_{\infty}^s, \qquad \underline{G}_n=V_{h^{(\mathrm{max})}}^{-1}
\llVert \mathbf{K}\rrVert _{\infty
}^s, \qquad \overline{G}_n=V_{h^{(\mathrm{min})}}^{-1}
\llVert \mathbf{K}\rrVert _{\infty}^s,\\
  G_{j,n}(h_j)&=&
\frac{h_j^{(\mathrm{min})}}{h_j} V_{h^{(\mathrm{min})}}^{-1}\llVert
\mathbf{K}\rrVert
_{\infty}^s,
\\
\underline{G}_{j,n}&=&\frac{h_j^{(\mathrm{min})}}{h_j^{(\mathrm{max})}}
 V_{h^{(\mathrm{min})}}^{-1}
\llVert \mathbf{K}\rrVert _{\infty}^s,\qquad  j=\overline {1,s},\qquad
\varrho_n^{(s)} (\widehat{h},\bar{h} )=\max
_{j=\overline{1,s}}\bigl\llvert \ln(\widehat{h}_j)-\ln(
\bar{h}_j)\bigr\rrvert .
\end{eqnarray*}

Obviously, Assumption $1$(i) in Lepski \cite{lepski:upperfunctions} is
fulfilled. Using Assumption (\ref{eq:kernelasumption}) (see
Section~\ref{sec:kernel} of the present paper), we get $\operatorname{supp}(K)\subseteq
[-1/2,1/2]^s$ and
\[
\bigl\llvert K(x)-K(y)\bigr\rrvert \leq L_\mathbf{K}^{(s)}\max
_{j=\overline
{1,s}}\llvert x_j-y_j\rrvert \qquad
\forall x,y\in\bR^s,\qquad  L_\mathbf {K}^{(s)}:=s\llVert
\mathbf{K}\rrVert _{\infty}^{s-1}L_{\mathbf{K}}>0. %
\]

Thus, we easily check that, for any $h,h'\in\cH_n^{(s)}$ and any
$y\in
\bR^s$,
\begin{eqnarray*}
&&\bigl\llvert K_h(y-y_0)-K_{h'}(y-y_0)
\bigr\rrvert \\
&&\quad\leq \biggl[\frac{\llVert \mathbf{K}\rrVert _{\infty}^s}{V_h}\vee \frac
{\llVert \mathbf{K}\rrVert _{\infty}^s}{V_{h'}} \biggr] \biggl
\{\exp \bigl(s\varrho_n^{(s)}\bigl(h,h'\bigr)
\bigr)-1+\frac{L_\mathbf{K}^{(s)}}{\llVert
\mathbf{K}\rrVert _{\infty}^s} \bigl(\exp \bigl(\varrho _n^{(s)}
\bigl(h,h'\bigr) \bigr)-1 \bigr) \biggr\}. %
\end{eqnarray*}
It implies that Assumption $1$(ii) in Lepski \cite{lepski:upperfunctions} holds with
\[
D_0(z)=\exp{(sz)}-1+\frac{L_\mathbf{K}^{(s)}}{\llVert \mathbf
{K}\rrVert
_{\infty}^s}\times \bigl(\exp{(z)}-1
\bigr),\qquad  D_{s+1}\equiv0, L_{s+1}\equiv0. %
\]

Furthermore, Assumption $3$ in Lepski \cite{lepski:upperfunctions}
holds with $N=0$ and $R=1$ since $\mH_{k+1}^m=\mH_{s+1}= \{
y_0
\}$ and Assumption $2$ in Lepski \cite{lepski:upperfunctions} is not
needed since $n_1=n_2=n$.

(2) Thus, the application of the Theorem $1$ in Lepski
\cite{lepski:upperfunctions} is possible. Let us first compute the
constants which appear in its proof.
\begin{eqnarray*}
C_{N,R,m,k}&=&\sup_{\delta>\delta_*}\delta^{-2}s \biggl[1+
\ln{ \biggl(\frac
{9216(s+1)\delta^{2}}{[s^*(\delta)]^{2}} \biggr)} \biggr]_+ +\sup_{\delta
>\delta_*}
\delta^{-2}s \biggl[1+\ln{ \biggl(\frac{9216(s+1)\delta
}{[s^*(\delta)]} \biggr)}
\biggr]_+\\
&:=&C_s; %
\\
C_D&=&se^s+\frac{seL_\mathbf{K}}{\llVert \mathbf{K}\rrVert _{\infty}},\qquad
C_{D,b}=
\sqrt{2C_D}\vee \bigl[(2/3) (C_D\vee8e) \bigr], \\
 \lambda
_1&=&4\sqrt {2eC_D}, \lambda_2=(16/3)
(C_D\vee8e). %
\end{eqnarray*}

Next, we have to compute the quantities involved in the description of
$\cU_{\mathbf{r}}^{(u,q)}(n,\mathfrak{h})$.
\[
M_q(h)\leq C_{s,1}^{(q)} \biggl[1\vee\ln \biggl(
\frac
{V_{h^{(\mathrm{max})}}}{V_h} \biggr) \biggr],\qquad  C_{s,1}^{(q)}:=\bigl[144s\delta _*^{-2}+5q+3+36C_s\bigr]\vee1.
\]

Since $Y_i, i=\overline{1,n}$, are identically distributed, putting
$\mathfrak{h}
=(h,y_0)$, $n_1=n_2=n$ and $\mathbf{r}=0$, we have
\begin{eqnarray*}
F_{n,\mathbf{r}}(\mathfrak{h})&=&1\vee \biggl[\int_{\bR^s}{
\bigl\llvert K_h(y-y_0)\bigr\rrvert g(y)\,\rd y}
\biggr]:=G_h(y_0),\\
  F_n&=&\sup
_{h\in\cH
_n^{(s)}}G_h(y_0)\leq1\vee\mathbf{g}
\llVert \mathbf{K}\rrVert _{1}^s; %
\\
\cU_{\mathbf{r}}^{(u,q)}(n,\mathfrak{h}) &\leq&\cU^{(u,q)}(n,h,y_0),\qquad
c(\mathbf{K},s,q):=\bigl[(10C_D)\vee (48e)\bigr]C_{s,1}^{(q)}
\llVert \mathbf{K}\rrVert _{\infty}^s. %
\end{eqnarray*}
Here, we have used that $C_{s,1}^{(q)}\wedge\llVert \mathbf{K}\rrVert
_{\infty}^s\geq1$. Thus, we come to the inequality (\ref
{eq:empiricalupperbound1}) with $C_s^{(q)}(\mathbf{K},\mathbf
{g}):=c_q\llVert \mathbf{K}\rrVert _{\infty}^{sq} (1\vee\mathbf
{g}\llVert \mathbf{K}\rrVert _{1}^s )^{q/2}$,
$c_q=2^{7q/2+5}3^{q+4}\Gamma(q+1)(C_{D,b})^q$.

(3) If $n\geq3, nV_h\geq\ln(n), 1\leq u\leq q\ln
(n)$ and $M(h):=1\vee\ln (\frac{V_{h^{(\mathrm{max})}}}{V_h} )$, since
$1\leq G_h(y_0)\leq\llVert \mathbf{K}\rrVert _{\infty}^s$, one has
%
\begin{eqnarray}
\label{eq:empiricalupperbound2} (nV_h)^{-1} \bigl\{M(h)+2\ln \bigl(2+\ln
G_h(y_0) \bigr)+u \bigr\} &\leq&7(nV_h)^{-1}G_h(y_0)
\bigl\{M(h)+u \bigr\}
\nonumber
\\[-8pt]
\\[-8pt]
\nonumber
&\leq&7(1+q)\llVert \mathbf {K}\rrVert _{\infty}^s.
\end{eqnarray}

Put finally $\lambda_s^{(q)}[\mathbf{K}]:=c(\mathbf{K},s,q)\sqrt
{7}
\{\sqrt{7(1+q)\llVert \mathbf{K}\rrVert _{\infty}^s}+1 \}$. Since
$ [a_s^{(q)} ]^{-1}\geq1$, the assertion~(i) of Proposition~\ref
{prop:empiricalupperbound1} follows from (\ref
{eq:empiricalupperbound1}) and (\ref{eq:empiricalupperbound2}). Let us
now prove the assertions (ii) and (iii) of Proposition~\ref
{prop:empiricalupperbound1}.

(4) First, in view of the definition of $\mH
_s^{(q)}(n)$, we get the assertion (ii) from the assertion (i) of
Proposition~\ref{prop:empiricalupperbound1} since $u\leq q\ln(n)$ and
$ [1\vee\lambda_s^{(q)} ]\sqrt{(1+q)a_s^{(q)}}=1/2$.
Here, we
have used that if $\mathbf{K}$ satisfies the assumption (\ref
{eq:kernelasumption}), see Section~\ref{sec:kernel}, $\llvert \mathbf
{K}\rrvert $ satisfies it as well and, therefore, Proposition~\ref
{prop:empiricalupperbound1}(i) is applicable to the process $\overline
{\xi}_h^{(n)}(y_0)$.

Next, using the trivial inequality $\llvert x\vee a-x\vee b\rrvert \leq
\llvert a-b\rrvert $, $x, a, b\in\bR$, we easily check that
%
\begin{eqnarray}
\label{eq:empiricalupperbound21} G_h(y_0)&\leq& 2\widetilde{G}_h(y_0)+2
\sup_{h\in\mH
_s^{(q)}(n)} \biggl[\bigl\llvert \overline{\xi}_h^{(n)}(y_0)
\bigr\rrvert -\frac{1}{2}G_h(y_0) \biggr]_+\qquad
\forall h\in\mH_{a_s}^{(s)}(n).
\end{eqnarray}
Assertion (iii) of Proposition~\ref{prop:empiricalupperbound1} follows
from assertion (ii) and (\ref{eq:empiricalupperbound21}).

\subsection{Proof of Lemma \texorpdfstring{\protect\ref{lem:empiricalupperbound}}{1}}

Note first that, for any $(h,\cP)\in\overline{\mH}[ \overline{\mP}
]$, any $(\eta,\cP')\in\overline{\mH}[ \overline{\mP} ]$ and any
$I\cap I'\in\cP\circ\cP'$
\begin{eqnarray*}
h_{I\cap I'}\vee\eta_{I\cap I'}&\in&\bigcup
_{m=1}^{M_n(I)}\bigcup_{l=1}^{M_n(I')}
\mH_{m,l}^{(I\cap I')},\\
 \mH_{m,l}^{(I\cap I')}&:=& \biggl
\{h_{I\cap I'}\in\prod_{i\in I\cap I'} \biggl[
\frac{1}{n},\mathfrak {h}_i^{(I\cap
I', m, l)} \biggr]\dvt
nV_{h_I}\geq \bigl[a_{\llvert I\cap I'\rrvert }^{(2q)}
\bigr]^{-1}\ln(n) \biggr\}, %
\end{eqnarray*}
where $\mathfrak{h}_i^{(I\cap I', m, l)}:=(2^{m\vee l})^{\mz}
[\mathfrak{h}_i^{(I)}\vee
\mathfrak{h}_i^{(I')} ]$, $i\in I\cap I'$.

Set $f\in\bF_d [\mathbf{f},\overline{\mP}  ]$. To get the
assertions of Lemma~\ref{lem:empiricalupperbound}, we apply Proposition~\ref{prop:empiricalupperbound1} with $s=\llvert I\cap I'\rrvert $,
$g=f_{I\cap I'}$, $\mathbf{g}=\mathbf{f}$, $h_i^{(\mathrm{min})}(n)=\frac
{1}{n}$, $h_i^{(\mathrm{max})}(n)=\mathfrak{h}_i^{(I\cap I', m, l)}$, $\mH
_s^{(q)}(n)=\mH
_{m,l}^{(I\cap I')}$, $K_{h}=K_{h_{I\cap I'}}^{(I\cap I')}$,
$G_h(y_0)=G_{h_{I\cap I'}}(x_{0,I\cap I'})$, $\widetilde
{G}_h(y_0)=\widetilde{G}_{h_{I\cap I'}}(x_{0,I\cap I'})$, $\cU
_h^{(u)}(y_0)=\cU_{h_{I\cap I'}}^{(u)}(x_{0,I\cap I'})$, $\xi
_{h}^{(n)}(y_0)=\xi_{h_{I\cap I'}}^{(n)}(x_{0,I\cap I'})$.

Recall that $\overline{\mP}^*:= \{\cP\circ\cP^{\prime}\dvt \cP, \cP
^{\prime}\in\overline{\mP} \}$. In view of the definition of
$\overline{\mH}[ \overline{\mP} ]$, we easily check that
\[
\xi_n(x_0)\leq\sum_{\cP\circ\cP'\in\overline{\mP}^*}
\sum_{I\cap I'\in\cP
\circ\cP'}\sum_{m=1}^{M_n(I)}
\sum_{l=1}^{M_n(I')}\sup_{h_{I\cap
I'}\in
\mH_{m,l}^{(I\cap I')}}
\bigl[\bigl\llvert \xi_{h_{I\cap
I'}}^{(n)}(x_{0,I})\bigr
\rrvert -\lambda_{\llvert I\cap I'\rrvert }^{(2q)}\cU _{h_{I\cap I'}}^{(u)}(x_{0,I\cap I'})
\bigr]_+, %
\]
with $u=q [1\vee\ln (2^{m\wedge l}V_{\mathrm{max}}/\inf_{I\in\cP
}V_{\mathfrak{h}
_I^{(I)}} ) ]\in[1,2q\ln(n)]$, since $V_{\mathfrak
{h}_I^{(I)}}\geq\frac
{\ln(n)}{an}$ and $M_n(I)\leq\log_2(n)$, $\forall I\in\cI_d$.

Therefore, it follows from the assertion (i) of Proposition~\ref
{prop:empiricalupperbound1}, since $V_{\mathfrak{h}_{I\cap
I'}^{(I)}\vee\mathfrak{h}
_{I\cap I'}^{(I')}}\geq\inf_{I\in\cP}V_{\mathfrak{h}_I^{(I)}}$,
\begin{eqnarray*}
&& \Bigl(\bE_f^{(n)} \Bigl\{\sup_{h_{I\cap I'}\in\mH_{m,l}^{(I\cap
I')}}
\bigl[\bigl\llvert \xi_{h_{I\cap I'}}^{(n)}(x_{0,I})\bigr
\rrvert -\lambda _{\llvert I\cap I'\rrvert }^{(2q)}\cU_{h_{I\cap I'}}^{(u)}(x_{0,I\cap I'})
\bigr]_+ \Bigr\}^{2q} \Bigr)^{{1}/{(2q)}}
\\
&&\quad\leq \bigl\{C_{\llvert I\cap I'\rrvert }^{(2q)}(\mathbf{K},\mathbf {g}) \bigr
\}^{{1}/{(2q)}} [nV_{\mathrm{max}} ]^{-{1}/{2}} \bigl(2^{\mz\llvert I\cap I'\rrvert /2}
\bigr)^{-m\vee l} \bigl(2^{1/2} \bigr)^{-m\wedge l};
\\
&&\bigl(\bE_f^{(n)}\bigl\llvert \xi_n(x_0)
\bigr\rrvert ^{2q} \bigr)^{
{1}/{(2q)}}\leq \mathbf{c}_1
[nV_{\mathrm{max}} ]^{-{1}/{2}},\\
 \mathbf {c}_1&:=&\sum
_{\cP\in\overline{\mP}^*}\sum_{I\in\cP} \bigl
\{C_{\llvert I\rrvert }^{(2q)}(\mathbf{K},\mathbf{f}) \bigr\}^
{{1}/{(2q)}}
\biggl[\frac
{2^{[(\mz|I|)\wedge1]/2}}{2^{[(\mz|I|)\wedge1]/2}-1} \biggr]. %
\end{eqnarray*}

Similarly, applying Proposition~\ref{prop:empiricalupperbound1}(iii)
and using the trivial inequality $[\sup_ix_i-\sup_iy_i]_+\leq\sup_i[x_i-y_i]_+$, we obtain the assertion (ii) of Lemma~\ref
{lem:empiricalupperbound} with $\mathbf{c}_2:=2\mathbf{c}_1$.

Next, it is easily seen that
\begin{eqnarray*}
\overline{G}_{n}(x_{0})&\leq&2 \Biggl(\sum
_{\cP\circ\cP'\in
\overline{\mP
}^*}\sum_{I\cap I'\in\cP\circ\cP'}\sum
_{m=1}^{M_n(I)}\sum_{l=1}^{M_n(I')}
\sup_{h_{I\cap I'}\in\mH_{m,l}^{(I\cap I')}} \biggl[\bigl\llvert \overline{\xi}_{h_{I\cap I'}}^{(n)}(x_{0,I})
\bigr\rrvert -\frac
{1}{2}G_{h_{I\cap I'}}(x_{0,I}) \biggr]_+
\Biggr)\\
&&{}+3G(x_0), %
\end{eqnarray*}
and that
\[
\bigl(\bE_f^{(n)}\bigl\llvert \overline{
\mathbf{f}}_n(x_0)\bigr\rrvert ^{2q}
\bigr)^{1/2q} \leq12\lambda d^3 2^{d^2} \bigl[ \bigl(
\bE_f^{(n)}\bigl\llvert \overline {G}_{n}(x_{0})
\bigr\rrvert ^{2qd^2} \bigr)^{{1}/{(2qd^2)}} + \bigl(1\vee\mathbf{f}
\llVert \mathbf{K}\rrVert _1^d \bigr)
\bigr]^{d^2}. %
\]
Thus, we get assertion (iii) of Lemma~\ref{lem:empiricalupperbound}
from assertion (ii) of Proposition~\ref{prop:empiricalupperbound1} with
\[
\mathbf{c}_3:=12\lambda d^3 \biggl[4 \biggl(\sum
_{\cP\in\overline
{\mP
}^*}\sum_{I\in\cP\circ\cP'} \bigl
\{C_{\llvert I\rrvert }^{(2qd^2)}(\mathbf {K},\mathbf{f}) \bigr\}
^{{1}/{(2qd^2)}}
\biggl[\frac{2^{[(\mz
|I|)\wedge
1]/2}}{2^{[(\mz|I|)\wedge1]/2}-1} \biggr] \biggr) +8 \bigl(1\vee\mathbf{f}\llVert
\mathbf{K}\rrVert _1^d \bigr) \biggr]^{d^2}.
\]

Similarly, we obtain assertion (iv) of Lemma~\ref
{lem:empiricalupperbound} with
\[
\mathbf{c}_4:=2 \biggl(\sum_{\cP\in\overline{\mP}^*}\sum
_{I\in
\cP\circ\cP
'} \bigl\{C_{\llvert I\rrvert }^{(2q)}(
\mathbf{K},\mathbf{f}) \bigr\} ^{{1}/{(2q)}} \biggl[\frac{2^{[(\mz|I|)\wedge1]/2}}{2^{[(\mz|I|)\wedge
1]/2}-1} \biggr]
\biggr) +3 \bigl(1\vee\mathbf{f}\llVert \mathbf{K}\rrVert _1^d
\bigr). %
\]
This completes the proof of Lemma~\ref{lem:empiricalupperbound}.

\subsection{Proof of Lemma \texorpdfstring{\protect\ref{lem:biasupperbound}}{3}} The proof
of this lemma is based on the embedding theorem for anisotropic
Nikolskii classes; see, for example, Theorem~6.9 in Nikolskii \cite{nikolskii}.

Let $\cP'\in\overline{\mP}$ and $I\in\cP\circ\cP'$ be fixed.
Set $\kappa
(I):=1-\sum_{k\in I} (\beta_kp_k )^{-1}$ and $\beta
_i(I):=\kappa(I)\beta_i\kappa_i^{-1}(I)$, where $\kappa
_i(I):=1-\sum_{k\in I} (p_k^{-1}-p_i^{-1} )\beta_k^{-1}$, $i\in I$. Since
$\kappa(I)>0$ there exists $c_I:=c_I(\mathbf{K},\llvert I\rrvert ,p_I,l)>0$ such that
\[
\bN_{p_I,\llvert I\rrvert }(\beta_I,L_I)\subseteq
\bN_{\infty,\llvert I\rrvert }\bigl(\beta(I),c_IL_I\bigr).
\]

Introduce the family of $\llvert I\rrvert \times\llvert I\rrvert $ matrices
$E_j:=(e_1,\ldots,e_j,0,\ldots,0), j=\overline{1,\llvert I\rrvert }$,
and $E_0$ is zero matrix. For any $(h,\eta)\in(0,1]^d\times[0,1]^d$,
using a telescopic sum and the triangle inequality, we get
\begin{eqnarray*}
&&\bigl\llvert \cB_{h_I,\eta_I}(x_{0,I})\bigr\rrvert \leq\sum
_{j=1}^{\llvert I\rrvert } \biggl|\int K^{(I)}(u)
\bigl[f_I \bigl(x_{0,I}+\eta_I
u+(h_I\vee\eta _I-\eta _I)E_ju
\bigr)\\
&&\hspace*{91pt}{}-f_I \bigl(x_{0,I}+\eta_I
u+(h_I\vee\eta_I-\eta_I)E_{j-1}u
\bigr)\bigr]\,\rd u\biggr |.
\end{eqnarray*}

For $j=1,\ldots,\llvert I\rrvert $ put
\begin{eqnarray*}
\cB_{h_I,\eta_I,j}(x_{0,I})&:=&\int_{\bR}
\mathbf{K}(u_j)\bigl[f_I \bigl(x_{0,I}+
\eta_I u+(h_I\vee\eta_I-
\eta_I)E_ju \bigr)
 \\
 &&\hspace*{41pt}{}-f_I \bigl(x_{0,I}+\eta_I
u+(h_I\vee\eta_I-\eta_I)E_{j-1}u
\bigr)\bigr]\,\rd u_j.
\end{eqnarray*}
If $\eta_j\geq h_j$, then $\cB_{h_I,\eta_I,j}(x_{0,I})=0$, if not we
put $[u]^j:=u-u_je_j, u\in\bR^{\llvert I\rrvert }$, and we have
\begin{eqnarray*}
\cB_{h_I,\eta_I,j}(x_{0,I})&=&\int_{\bR}
\mathbf{K}(u_j) \bigl[f_I \bigl(x_{0,I}+
\eta_I u+(h_I\vee\eta_I-
\eta_I)E_ju \bigr)
\\
&&\hspace*{41pt}{} -f_I \bigl(x_{0,I}+ [\eta_I u
]^j+(h_I\vee\eta_I-\eta_I)E_{j-1}u
\bigr) \bigr]\,\rd u_j
\\
&&{}+\int_{\bR} \mathbf{K}(u_j)
\bigl[f_I \bigl(x_{0,I}+ [\eta_I u
]^j+(h_I\vee\eta_I-\eta_I)E_{j-1}u
\bigr) \\
&&\hspace*{52pt}{}-f_I \bigl(x_{0,I}+\eta_I
u+(h_I\vee\eta_I-\eta_I)E_{j-1}u
\bigr) \bigr]\,\rd u_j.
\end{eqnarray*}

Thus, in view of the triangle inequality,
\begin{eqnarray*}
\bigl\llvert \cB_{h_I,\eta_I}(x_{0,I})\bigr\rrvert &\leq&2\sum
_{i\in
I}c_IL_ih_i^{\beta_i(I)}
\int_{\bR^{\llvert I\rrvert }}{\bigl\llvert K^{(I)}(u)\bigr\rrvert
\llvert u_i\rrvert ^{\beta_i(I)}}\,\rd u \leq\mathbf{c}\sum
_{i\in I}L_ih_i^{\beta_i(I)},
\\
\mathbf{c}&:=&\mathbf{c}(\mathbf{K},d,p,l,\cP)=2\llVert \mathbf {K}\rrVert
_1^{d}\sup_{\cP'\in\overline{\mP}}\sup
_{I\in\cP\circ\cP
'}c_I\bigl(\mathbf {K},\llvert I\rrvert
,p_I,l\bigr).
\end{eqnarray*}
Here, we have used Taylor expansions of $f\in\bN_{\infty,\llvert I\rrvert }(\beta(I),c_IL_I)$, the product structure of $K^{(I)}$, the Fubini
theorem that $\beta(I)\in(0,l]^d$ and (\ref{eq:kernelorthogonality});
see Section~\ref{sec:minimaxresults}. We have also used that $\mathbf
{K}$ is compactly supported on $[-1/2,1/2]$ and that $\llVert \mathbf
{K}\rrVert _1\geq1$.

\subsection{Proof of Lemma \texorpdfstring{\protect\ref{lem:adaptivelowerbound}}{4}} Put
$T_n:=a_n  |\widetilde{f}_n(x_0)-f_0(x_0) |$ and
\[
\cR_n^{(q)} [a_n,b_n,
\widetilde{f},f ]:=\sup_{f\in
N_1}\bE _f^{(n)}
\bigl\{a_n \bigl|\widetilde{f}_n(x_0)-f(x_0)
\bigr| \bigr\}^q +\sup_{f\in N_2} \bE_f^{(n)}
\bigl\{b_n \bigl|\widetilde {f}_n(x_0)-f(x_0)
\bigr| \bigr\}^q. %
\]
It is easily seen that $\cR_n^{(q)} [a_n,b_n,\widetilde
{f},f
]\geq\cR_n^{(1)} [a_n,b_n,\widetilde{f},f ]$ and that
\begin{eqnarray*}
\cR_n^{(1)} [a_n,b_n,
\widetilde{f},f ]\geq\bE _{f_1}^{(n)} \bigl\{ |T_n-1 |
\bigr\}+\frac{b_n}{a_n}\bE_{f_0}^{(n)} \{T_n \}.
\end{eqnarray*}
Here, we have used the triangle inequality and the assumption $a_n
|f_1(x_0)-f_0(x_0) |=1$.

Put also $c_n:=\frac{b_n}{a_n}$ and $Z_n:=\frac{\rd\bP
_{f_1}^{(n)}}{\rd
\bP_{f_0}^{(n)}} (X^{(n)} )$. We obtain
\begin{eqnarray*}
\cR_n^{(1)} [a_n,b_n,
\widetilde{f},f ]\geq\bE_{f_0}^{(n)} \{c_n\wedge
Z_n \}\geq\tfrac{1}{2} \Bigl[c_n+1-\sqrt{
\bE _{f_0}^{(n)} \{c_n-Z_n
\}^2} \Bigr].
\end{eqnarray*}
Here, we have used the trivial equality $a\wedge b=\frac{1}{2} \{
a+b-|a-b| \}$, that $\bE_{f_0}^{(n)}  \{Z_n  \}=1$ and the
Cauchy--Schwarz inequality. Using the third assumption, we also have
$\bE_{f_0}^{(n)}  \{c_n-Z_n \}^2\leq c_n^2-c_n$.
Finally, for $n$ large enough,
\[
\inf_{\widetilde{f}}\cR_n^{(q)}
[a_n,b_n,\widetilde{f},f ]\geq \frac{1}{2}
\Bigl[c_n+1-\sqrt{c_n^2-c_n}
\Bigr]\geq\frac{1}{2}. %
\]
%
\end{appendix}

\section*{Acknowledgments} The author is grateful to O. Lepski and the
anonymous referees for their very useful remarks and suggestions.

%

%




\printhistory
\end{document}